\documentclass[12pt]{article}
\usepackage{amsmath,amssymb}

\newcommand{\ncm}{\newcommand}     
\renewcommand{\Re}{{\rm Re}}
 \ncm{\Inn}{\mbox{\rm Inn}} 
\ncm{\Ap}{\mbox{$\overline{\rm Inn}$}} \ncm{\Ext}{\mbox{\rm Ext}} 
\ncm{\Ex}{\mbox{\rm Ex}} \ncm{\OExt}{\mbox{\rm OrderExt}} 
\ncm{\AI}{\mbox{\rm AInn}} \ncm{\HI}{\mbox{\rm HInn($A$)}} 
\ncm{\Aut}{\mbox{\rm Aut}} \ncm{\Mal}{\mbox{$M_{\alpha}$}} 
\ncm{\Aff}{\mbox{${\rm Aff}$}} \ncm{\id}{\mbox{\rm id}} 
\ncm{\Ker}{\mbox{\rm Ker}} \ncm{\BE}{\begin{eqnarray*}} 
\ncm{\EE}{\end{eqnarray*}} \ncm{\lra}{\mbox{$\longrightarrow$}} 
\ncm{\Hom}{\mbox{\rm Hom}} \ncm{\calU}{{\cal U}} \ncm{\el}{\ell} 
 \ncm{\ad}{{\rm ad}} 
 \ncm{\diag}{{\rm diag}}    
 \ncm{\rank}{{\rm rank}}       
 \ncm{\Tr}{{\rm Tr}}  
\ncm{\Alg}{\mbox{\rm Alg}} \ncm{\Conv}{\mbox{\rm Conv}} 
\ncm{\D}{{\cal D}} 

\ncm{\cstar}{C$^{*}$-algebra} \ncm{\cstars}{C$^{*}$-algebras} 
\ncm{\subcstar}{C$^{*}$-subalgebra}  
\ncm{\ra}{\mbox{$\rightarrow$}} \ncm{\la}{\mbox{$\leftarrow$}} 
\ncm{\hra}{\hookrightarrow} \ncm{\da}{\mbox{$\downarrow$}} 
\ncm{\se}{\mbox{$\searrow$}} \ncm{\al}{\mbox{$\alpha $}} 
\ncm{\del}{\mbox{$\delta$}} \ncm{\supp}{\mbox{\rm supp}} 
\ncm{\Ad}{\mbox{\rm Ad}} \ncm{\CAR}{\mbox{$M_{2^{\infty}}$}} 
\ncm{\ep}{\mbox{$\epsilon > 0$}} \ncm{\ol}{\overline} 
\ncm{\Mninf}{\mbox{$M_{n^{\infty}}$}} \ncm{\MR}{M. R\o{}rdam} 
\ncm{\Ran}{\mbox{\rm Ran}} 
\ncm{\vo}{}
\ncm{\ch}{}
\ncm{\CMP}{Comm. Math. Phys.} \ncm{\add}{} 
\ncm{\tilsig}{\tilde{\sigma}} \ncm{\dist}{{\rm 
dist}}\ncm{\eps}{\epsilon}  \ncm{\calL}{{\mathcal{L}}}   
 \ncm{\E}{{\mathcal{E}} }    \ncm{\M}{{\mathcal{M}} }
 \ncm{\calO}{{\mathcal{O}}} \ncm{\F}{{\mathcal{F}}} \ncm{\G}{{\mathcal{G}} }  
\ncm{\calH}{{\mathcal{H}}}   \ncm{\calC}{{\mathcal{C}}} 
\ncm{\lan}{{\langle}}\ncm{\ran}{{\rangle}}   
\ncm{\calK}{{\mathcal{K}}}    \ncm{\Spec}{{\rm Spec}}  
\ncm{\calP}{{\mathcal{P}}} \ncm{\Hil}{{\mathcal{H}}} 
\ncm{\U}{{\mathcal{U}}}                \ncm{\A}{\mathcal{A}}    
\ncm{\Nr}{\mathcal{N}}     

\newtheorem{theo}{Theorem}[section]

\newtheorem{lem}[theo]{Lemma}
\newtheorem{prop}[theo]{Proposition}
\newtheorem{remark}[theo]{Remark}
\newtheorem{definition}[theo]{Definition}
\newtheorem{example}[theo]{Example}
\newtheorem{property}[theo]{Property}

\newenvironment{rem}{\begin{remark} \rm}{\end{remark}}
\newenvironment{pf}{{\it Proof.}}{\hfill$\square$\vspace{3mm}}

\ncm{\R}{\mbox{\bf R}} \ncm{\Z}{\mbox{\bf Z}} \ncm{\T}{\mbox{\bf 
T}} \ncm{\TT}{\T$^{2}$} \ncm{\N}{\mbox{\bf N}} \ncm{\C}{\mbox{\bf 
C}}

\title{Approximate AF flows}
\bigskip
\author{Akitaka Kishimoto\\  
Department of Mathematics, Hokkaido University \\
 Sapporo, 060-0810 Japan\\ 
\\
{\em Dedicated to George A. Elliott on the occasion of his 
sixtieth birthday}} 
\date{}
\begin{document}
\maketitle 

\begin{abstract}
When $\alpha$ is a flow on a unital AF algebra $A$ such that there 
is an increasing sequence $(A_n)$ of finite-dimensional 
$\alpha$-invariant C$^*$-subalgebras of $A$ with dense union, we 
call $\alpha$ an AF flow. We show that an {\em approximate} AF 
flow is a cocycle perturbation of an AF flow. 
\end{abstract}                  
                    
\section{Introduction}
Let $\alpha$ be a flow on a unital \cstar\ $A$, i.e., $t\mapsto 
\alpha_t$ is a group homomorphism of the real line $\R$ into the 
automorphism group of $A$ such that $t\mapsto \alpha_t(x)$ is 
continuous for each $x\in A$. If $u$ is a continuous map of $\R$ 
into the unitary group $\U(A)$ of $A$ such that 
$u_s\alpha_s(u_t)=u_{s+t},\ s,t\in\R$, then $u$ is called an 
$\alpha$-cocycle. In this case $t\mapsto \Ad\,u_t\alpha_t$ is a 
flow again and is called a cocycle perturbation of $\alpha$. 

We denote by $\delta_\alpha$ the generator of $\alpha$, which is a 
closed derivation defined on a dense $*$-subalgebra 
$\D(\delta_\alpha)$ \cite{BR1,Sak}. If the $\alpha$-cocycle $u$ is 
differentiable  and $h=-idu_t/dt|_{t=0}\in A_{sa}$, then the 
generator of the flow $t\mapsto \Ad\,u_t\alpha_t$ is given by 
$\delta_\alpha+\ad\,ih$. If the $\alpha$-cocycle $u$ is a 
coboundary, i.e., $u_t=w\alpha_t(w^*)$ for some $w\in\U(A)$, then 
the generator of the flow $t\mapsto \Ad\,u_t\alpha_t$ is given by 
$\Ad\,w\circ\delta_\alpha\circ\Ad\,w^*$ on 
$\Ad\,w(\D(\delta_\alpha))$. In general the $\alpha$-cocycle is 
given as a combination of the above two types \cite{Kis00}.  

If $A$ is an AF algebra \cite{Br} and has an increasing sequence 
$(A_n)$ of finite-dimensional C$^*$-subalgebras of $A$ such that 
$\bigcup_nA_n$ is dense in $A$ and $\alpha_t(A_n)=A_n$ for all 
$t\in\R$ and all $n$, then $\alpha$ is called an {\em AF flow}. 
Note that in this case there is a self-adjoint element $h_n\in 
A_n$ such that $\alpha_t|A_n=\Ad\,e^{ith_n}|A_n$. It also follows 
that $\D(\delta_\alpha)\supset A_n$, 
$\delta_\alpha|A_n=\ad\,ih_n|A_n$, and $\bigcup_nA_n$ is a core 
for $\delta_\alpha$. Since $[h_m,h_n]=0$, the generator of this 
type is called {\em commutative} and is studied in Sakai's book 
\cite{Sak} (see also \cite{Kis00,BK,De}). In particular $\alpha$ 
is approximately inner in the sense that 
$\lim_n\Ad\,e^{ith_n}(x)=\alpha_t(x)$ uniformly in $t$ on every 
bounded subset of $\R$ for all $x\in A$.             

AF flows (on a UHF algebra) appear as time-flows for classical 
lattice models in physics and look manageable for analysis (e.g., 
the KMS states have explicit expressions \cite{Sak}). There are 
time-flows for quantum lattice models which are not obviously 
cocycle perturbations of AF flows, but we still lack a rigorous 
proof to that effect though we know that there are flows which are 
not cocycle perturbations of AF flows. Our main concern is to 
distinguish the class of cocycle perturbations of AF flows among 
the flows which occur in physical models and thus to understand 
the flows beyond this class better. In this note we give a 
characterization of this class. 

When $B$ and $C$ are C$^*$-subalgebras of $A$, we write 
$B\stackrel{\delta}{\subset}C$ if for any $x\in B$ there is $y\in 
C$ such that $\|x-y\|\leq \delta \|x\|$. We define the distance of 
$B$ and $C$ by
 $$
 \dist(B,C)=\inf\{\delta>0\ |\ B\stackrel{\delta}{\subset}C,\ 
 C\stackrel{\delta}{\subset}B\}.
 $$                              

If $\alpha$ is an AF flow, then a cocycle perturbation $\alpha'$ 
of $\alpha$ may not be an AF flow but an {\em approximate AF flow} 
in the sense that $\sup_{t\in[0,1]}\dist(\alpha'_t(A_n),A_n)\ra0$, 
where the sequence $(A_n)$ is chosen for $\alpha$ as above. Our 
main purpose is to show the converse:

\begin{theo}  \label{M}
Let $\alpha$ be a flow on a unital AF algebra $A$. Then the 
following conditions are equivalent: 
 \begin{enumerate}
 \item $\alpha$ is a cocycle perturbation of an AF flow.
 \item $\alpha$ is an approximate AF flow, i.e., there 
 is an increasing sequence $(A_n)$ of finite-dimensional
 C$^*$-subalgebras of $A$ such that $\bigcup_nA_n$ is dense in $A$ and 
 $$
 \sup_{t\in [0,1]}\dist(\alpha_t(A_n),A_n)\ra0
 $$ 
 as $n\ra\infty$.
 \end{enumerate}
\end{theo}  

As noted above, the former implies the latter. In the rest of this 
note we shall prove that the latter implies the former. Since the 
latter condition is preserved under cocycle perturbations, it 
suffices to show: 

\begin{lem} \label{M1}    
For any $\eps>0$ there exists a $\delta>0$ satisfying the 
following condition: If $\alpha$ is an approximate AF flow, i.e., 
there is an increasing sequence $(A_n)$ of finite-dimensional 
C$^*$-subalgebras of $A$ such that $\bigcup_nA_n$ is dense in $A$ 
and 
 $$
 \delta_n\equiv\sup_{t\in [0,1]}\dist(\alpha_t(A_n),A_n)\ra0
 $$ 
as $n\ra\infty$ with $\delta_1\leq \delta$, there is an 
$\alpha$-cocycle $u$ such that $\|u_t-1\|<\eps$ for $t\in [0,1]$ 
and $\Ad\,u_t\alpha_t(A_1)=A_1$.   

Moreover if $A_0$ is a C$^*$-subalgebra of $A_1$ such that 
$\alpha_t(A_0)=A_0$, then $u$ can be chosen from $A\cap A_0'$. 
\end{lem} 
          
\begin{rem}
For a single automorphism $\alpha$ of a unital AF algebra $A$ the 
following conditions are equivalent:
 \begin{enumerate}
 \item For any $\eps>0$ there is a $u\in\U(A)$ and an increasing sequence 
 $(A_n)$ of finite-dimensional C$^*$-subalgebras of $A$ such that $\bigcup_nA_n$
 is dense in $A$ and $\Ad\,u\,\alpha(A_n)=A_n$ for all $n$.
 \item There is an increasing sequence $(A_n)$ of finite-dimensional 
 C$^*$-subalgebras of $A$ such that $\bigcup_nA_n$
 is dense in $A$ and  
 $$
 \dist(\alpha(A_n),A_n)\ra0
 $$                        
 as $n\ra\infty$.
 \end{enumerate}
This follows from a deep result of Christensen's (see 5.3 of
\cite{Ch}), which will be also used in the proof of our main 
result. 

There is a totally different sufficient condition for $\alpha$ 
that implies the above condition 1, i.e., which says that $\alpha$ 
has the Rohlin property and $\alpha_*$ is the identity on the 
dimension group $K_0(A)$. See \cite{Voi,Kis95,EK}. It is certainly 
desirable to find a sufficient condition like this for AF flows. 
\end{rem} 

\begin{rem}
Theorem \ref{M} might hold for non-unital AF algebras, where the 
$\alpha$-cocycle should be understood as a function into the 
multiplier algebra with continuity for the strict topology. But in 
this case it follows, by the same proof, that 2 implies 1, but it 
is not obvious how to prove that 1 implies 2. 
\end{rem}  
                                                                
The main theorem \ref{M} is what we should have settled sooner or 
later after singling out AF flows, but is not likely to be useful 
to distinguish the class of cocycle perturbations of AF flows.  
However there are some other attempts to characterize this class. 
In section 2 we will briefly survey them. In section 3 we will 
present a key idea for proving the main theorem in the setting of 
matrix algebras, which is a degenerate case of Prop.~\ref{M3} for 
the UHF algebra (see below). Letting $M_N$ to be the \cstar\ of 
$N\times N$ matrices, what will be shown is as follows:

\begin{prop}   \label{A}
For any $\eps>0$ there exists a $\delta>0$ satisfying the 
following condition: Let $A=M_N$ for some $N\in \N$, $B$ a unital 
$*$-subalgebra of $A$ with $B\cong M_K$ for some $K$ dividing $N$, 
and $H\in A_{sa}$ such that 
 $$
 e^{itH}Be^{-itH}\stackrel{\delta}{\subset} B
 $$
for $t\in [0,1]$. Then there exist a $u\in \U(A)$ and $h\in 
A_{sa}$ such that $\|u-1\|<\eps$, $\|h\|<\eps$, and 
 $$
 u^*e^{it(H+h)}uBu^*e^{-it(H+h)}u=B,\ \ t\in \R.
 $$
\end{prop}  

In the above situation let $\alpha_t=\Ad\,e^{itH}$, which defines 
a flow on $A=M_N$, and let $v_t=u^*e^{it(H+h)}ue^{-itH}$. Then 
note that $\Ad\,v_t\alpha_t(B)=B$ and that $v:t\mapsto v_t$ is an 
$\alpha$-cocycle such that $\sup_{t\in[0,1]}\|v_t-1\|<3\eps$. 

In section 4 we will give a version \ref{M3} of the main theorem 
for UHF flows \cite{Kis01a}, where $\alpha$ is a UHF flow on $A$ 
if $A$ is a UHF algebra and there is a sequence $(A_n)$ of unital 
matrix C$^*$-subalgebras of $A$ with dense union and 
$\alpha_t(A_n)=A_n,\ n\in\N,\ t\in\R$. Some of the technical 
results hold for more general unital \cstars. 

To go from UHF flows to AF flows, we will need a version of the 
following result, which is valid for any \cstars.           

\begin{prop}\label{Inv}
For any $\eps>0$ there exists a $\delta>0$ satisfying the 
following condition: Let $A$ be a unital \cstar\ and $\alpha$ be a 
flow on $A$. If $D$ is a unital finite-dimensional 
C$^*$-subalgebra of $A$ such that 
 $$
 \sup_{|t|\leq1}\|(\alpha_t-\id)|D\|<\delta,
 $$
then there is an $\alpha$-cocycle $u$ such that 
$\Ad\,u_t\alpha_t|D=\id$ and
 $$
 \max_{|t|\leq1}\|u_t-1\|<\eps.
 $$
\end{prop}  

What is important here is that $\delta$ does not depend on $D$ (or 
the size of $D$). In section 5 we will prove the above proposition 
and then derive the main result.

\section{AF flows}

For a \cstar\ $A$ we denote by $\ell^\infty(A)$ the \cstar\ of 
bounded sequences in $A$ and by $c_0(A)$ the ideal of 
$\ell^\infty(A)$ consisting of $x=(x_n)$ for which 
$\lim_{n\ra\infty}\|x_n\|=0$ and let 
$A^\infty=\ell^\infty(A)/c_0(A)$.  We embed $A$ into $A^\infty$ by 
regarding each $x\in A$ as the constant sequence $(x,x,\ldots)$. 
Given a flow $\alpha$ on $A$ we denote by $\ell^\infty_\alpha(A)$ 
the C$^*$-subalgebra of $x=(x_n)\in \ell^\infty(A)$ for which 
$t\mapsto \alpha_t(x)=(\alpha_t(x_n))$ is norm-continuous and 
define $A^\infty_\alpha$ as its image in $A^\infty$. We naturally 
have the flow $\ol{\alpha}$ on $A_\alpha^\infty$ induced by 
$\alpha$. We will also denote by $\alpha$ the restriction of 
$\ol{\alpha}$ to $A'\cap A_\alpha^\infty$.                                   

The following properties shared by AF flows and their cocycle 
perturbations could be used to distinguish them from other flows 
\cite{Kis01,BK}; similar properties are also considered for Rohlin 
flows \cite{Kis03,Kisp}.

\begin{prop}
Let $A$ be a unital AF algebra and let $\alpha$ be a cocycle 
perturbation of an AF flow on $A$. Then $(A^\infty_\alpha\cap 
A')^\alpha$ has real rank zero and has trivial $K_1$. Moreover 
$(A_\alpha^\infty)^\alpha$ has real rank zero and trivial $K_1$. 
\end{prop}
\begin{pf}   
The latter part is shown in 3.8 of \cite{BK} and 3.6 and 4.1 of 
\cite{Kis01}. The first part also follows similarly; but we will 
indicate how to prove it. 
                               
Apparently we may suppose that $\alpha$ is an AF flow. Hence we 
suppose that there is an increasing sequence $(A_n)$ of 
$\alpha$-invariant finite-dimensional $C^*$-subalgebras of $A$ 
with dense union. 
             
Let $b^*=b\in (A_\alpha^\infty\cap A')^\alpha$. Then there is a 
bounded sequence $(b_n)$ in $A_{sa}$ such that $b\sim (b_n)$ 
(i.e., $b=(b_n)+c_0(A)$). We may suppose that 
$\|\delta_\alpha(b_n)\|\ra0$ as $n\ra\infty$ and that there are 
increasing sequences $(k_n)$ and $(\ell_n)$ in $\N$ such that 
$k_n<\ell_n$, $k_n\ra\infty$, and $b_n\in B_n\equiv A_{\ell_n}\cap 
A_{k_n}'$. (The latter follows because $\bigcup_nA_n$ is a core 
for $\delta_\alpha$.) Since $B_n$ is $\alpha$-invariant and 
finite-dimensional, there is a $h_n^*=h_n\in B_n$ such that 
$\delta_\alpha|B_n=\ad\,h_n|B_n$. Since $\|[h_n,b_n]\|\ra0$ and 
$h_n,b_n\in (B_n)_{sa}$, we get $h_n',b_n'\in B_n$ such that 
$\|h_n-h_n'\|\ra0$, $\|b_n-b_n'\|\ra0$, and $[h_n',b_n']=0$ (see 
3.1 of \cite{BK}, which is an improvement of Lin's result 
\cite{Lin}).  

Let $\eps>0$ and let $F$ be a finite subset of the spectrum 
$\Spec(b)$ of $b$ such that any $\lambda\in \Spec(b)$ has $p\in F$ 
such that $|\lambda-p|<\eps$. Then we find a $b_n''\in (B_n)_{sa}$ 
such that $b_n''$ is a function of $b_n'$, 
$\limsup_n\|b_n'-b_n''\|<\eps$, and $\Spec(b_n'')\subset F$. Then 
$(b_n'')$ defines a self-adjoint element $c\in 
(A_\alpha^\infty\cap A')^\alpha$ such that $\|c-b\|<\eps$ and 
$\Spec(c)\subset F$, which is finite. This concludes the proof 
that $(A_\alpha^\infty\cap A')^\alpha$ has real rank zero. 

Let $u$ be a unitary in $(A_\alpha^\infty\cap A')^\alpha$. Then as 
before we may suppose that there is a sequence $(u_n)$ in $\U(A)$ 
and  increasing sequences $(\ell_n)$ and $(k_n)$ in $\N$ such that 
$k_n<\ell_n$, $k_n\ra\infty$, $u_n\in A_{\ell_n}\cap A_{k_n}'$, 
and $\|\delta_\alpha(u_n)\|\ra0$. There is an $h_n^*=h_n\in 
B_n\equiv A_{\ell_n}\cap A_{k_n}'$ such that 
$\delta_\alpha|B_n=\ad\,ih_n|B_n$. Then by using the condition 
that $\|[u_n,h_n]\|\ra0$, we apply 4.1 of \cite{Kis01}. 
\end{pf}                                    
          
There are some examples of approximately inner flows on an AF 
algebra without the above types of properties (see Section 3 of 
\cite{Kis01} and 3.11 of \cite{BK}, where only the properties for 
$(A_\alpha^\infty)^\alpha$ are explicitly mentioned). Those 
examples are of the following type. Let $C$ be a maximal abelian 
C$^*$-subalgebra (masa) of $A$ and choose a sequence $(h_n)$ in 
$C_{sa}$ such that the graph limit $\delta$ of $(\ad\,ih_n)$ is 
densely-defined and hence generates a flow \cite{BR1,Sak}. This is 
what we have as the examples and might be called a {\em quasi AF 
flow} (or a {\em commutative flow} following \cite{Sak}). Note 
that the domain $\D(\delta)$ of $\delta$ contains the masa $C$ 
(which is actually a Cartan masa in our examples); but depending 
on $(h_n)$ it may contain another masa as well. (We know of no 
example of a generator whose domain does not contain a masa.) 

\begin{rem} 
For a flow $\alpha$ on a unital simple AF algebra $A$ it is shown 
in \cite{Kis01} that $\alpha$ is a cocycle perturbation of an AF 
flow if and only if the domain $\D(\delta_\alpha)$ contains a 
canonical AF masa of $A$, where $C$ is a canonical AF masa if 
there is an increasing sequence $(A_n)$ of finite-dimensional 
C$^*$-subalgebras of $A$ with dense union such that $C$ is the 
closure of $\bigcup_nC\cap A_n\cap A_{n-1}'$ with $A_0=0$. 
\end{rem}

We note the following uniqueness result for canonical AF masas 
(cf. \cite{Ren}). 

\begin{prop}                                     
The canonical AF masas of an AF algebra are unique up to 
automorphism, i.e., if $A$ is an AF algebra and $C_1$ and $C_2$ 
are canonical AF masas of $A$, then there is an automorphism 
$\phi$ of $A$ such that $\phi(C_1)=C_2$. 
\end{prop} 
\begin{pf}
There are increasing sequences $(A_n)$ and $(B_n)$ of 
finite-dimensional C$^*$-subalgebras of $A$ such that 
$\bigcup_nA_n$ and $\bigcup_nB_n$ are dense in $A$ and $C_1$ 
(resp. $C_2$) is the closure $\bigcup_nC_1\cap A_n\cap A_{n-1}'$ 
(resp.$\bigcup_nC_2\cap B_n\cap B_{n-1}'$). By passing to 
subsequences, we find sequences $(u_n)$ and $(v_n)$ in $\U(A)$ 
such that $u_1=1$, $\|u_n-1\|<2^{-n}$, $\|v_n-1\|<2^{-n}$, 
$u_{n+1}\in \Ad(u_nu_{n-1}\cdots u_1)(A_n)'$, $v_{n+1}\in 
\Ad(v_n\cdots v_1)(B_n)'$, and 
 $$
 \Ad(u_n\cdots u_1)(A_n)\subset \Ad(v_n\cdots v_1)(B_n)\subset 
 \Ad(u_{n+1}u_n\cdots u_1)(A_{n+1}),
 $$
for all $n$. Let $u=\lim_nu_nu_{n-1}\cdots u_1$ and 
$v=\lim_nv_n\cdots v_1$. Then $u$ and $v$ are unitaries in $A$ (or 
$A+\C1$ if $A\not\ni 1$) such that $uA_nu^*\subset vB_nv^*\subset 
uA_{n+1}u^*$ for all $n$. Since $uC_1u^*$ and $vC_2v^*$ are also  
canonical AF masas, we may suppose that $A_n\subset B_n\subset 
A_{n+1}$. We choose maximal abelian C$^*$ subalgebras $D_n$ of 
$B_n\cap A_n'$ for $n=1,2,\ldots$ and $E_n$ of $A_{n+1}\cap B_n$ 
for $n=0,1,2,\ldots$ with $B_0=0$. Then the C$^*$-subalgebra $D$ 
generated by $D_n$ and $E_n$ for all $n$ is a canonical AF masa of 
$A$. Since the \cstar\ generated by $D_n$ and $E_n$ is isomorphic 
to $C_1\cap A_{n+1}\cap A_n'$ for $n=0,1,2,\ldots$ with $D_0=0$, 
there is a unitary $u_n\in A_{n+1}\cap A_n'+1$ such that 
$\Ad\,u_n(C_1\cap A_{n+1}\cap A_n')=C^*(D_n,E_n)$. Since the limit 
of $\Ad(u_0u_1u_2\cdots u_n)$ defines an automorphism of $A$, 
there is an automorphism $\phi_1$ of $A$ such that 
$\phi_1(C_1)=D$. Similarly, since the \cstar\ generated by $E_n$ 
and $D_{n+1}$ is isomorphic to $C_2\cap B_{n+1}\cap B_n'$, there 
is an automorphism $\phi_2$ such that $\phi_2(C_2)=D$. Thus 
$\phi_2^{-1}\phi_1(C_1)=C_2$, which concludes the proof. \end{pf}

\section{Matrix algebras} 

In this section we shall prove Proposition \ref{A}.

Let $\tau$ denote the unique tracial state of $A=M_N$, i.e., 
$\tau=(1/N)\Tr$, and define an inner product on $A$ by $\lan 
x,y\ran=\tau(y^*x)=\tau(xy^*), \ x,y\in A$. Equipped with this 
inner product, $A$ is an $N^2$-dimensional Hilbert space, which we 
will denote by $A_\tau$. We define a representation $\rho$ on 
$A_\tau$ of the tensor product $A\otimes A$ by $\rho(x\otimes 
y)\xi=x\xi y^t, \ \xi\in A_\tau$ and a representation $\pi$ of $A$ 
by $\pi(x)=\rho(x\otimes 1),\ x\in A$. Here $y^t$ denotes the 
tanspose of $y\in A=M_N$. Note that $\rho$ is irreducible and the 
state $\omega$ of $A\otimes A$ defined by $\omega(x\otimes y)=\lan 
\rho(x\otimes y)1,1\ran=\tau(xy^t)$ satisfies the condition that 
$\omega|B\otimes B^t$ is pure because $B\otimes B^t\cong M_{K^2}$ 
and the subspace $\rho(B\otimes B^t)1=B$ is $K^2$-dimensional.    

Let $U_t=\exp{ it(H\otimes1-1\otimes H^t)}=e^{it(H\otimes 
1)}e^{-it(1\otimes H^t)}$ and let $\gamma$ denote the flow 
$t\mapsto\Ad\,U_t$ on $A\otimes A$. Then $\gamma_s(x\otimes 
y^t)=\alpha_s(x)\otimes \alpha_s(y)^t$ for $x,y\in A$, where 
$\alpha$ is the flow on $A$ defined by 
$\alpha_t(x)=\Ad\,e^{itH}(x)$. We should note that $U$ has the 
following properties: $\rho(U_t)1=1$ and 
$\Ad\,\rho(U_t)\pi(x)=\pi\alpha_t(x),\ x\in A$.   
                                                     
If $\alpha_s(B)\stackrel{\delta}{\subset}B$ for $s\in[0,1]$ for a 
$\delta>0$ with $\delta<10^{-4}$, then Christensen \cite{Ch} shows 
that there is a unitary $v_s\in \U(A)$ such that 
$\|1-v_s\|<120\delta^{1/2}$ and 
 $$
 v_s\alpha_s(B)v_s^*=B.
 $$                                                   
Since $(v_s\otimes \ol{v}_s)(\alpha_s(B)\otimes 
\alpha_s(B)^t)(v_s\otimes \ol{v}_s)^*=B\otimes B^t$, we have that
 $$
 \gamma_s(B\otimes B^t)\stackrel{480\delta^{1/2}}{\subset}B\otimes 
 B^t,
 $$  
for $t\in [0,1]$.              

Hence the flow $\gamma$ on $A\otimes A$ satisfies the condition of 
Lemma \ref{A2} below for $(A,\alpha)$ with $B\otimes B^t$ and 
$\omega$ in place of $B$ and $\phi$ respectively if we start with 
the small enough $\delta>0$. Thus we get a $u\in \U(A\otimes A)$ 
and $h\in (A\otimes A)_{sa}$ such that $\|u-1\|\approx0$, 
$\|h\|\approx0$, and $\gamma'_t(B\otimes B^t)=B\otimes B^t,\ 
t\in\R$, where 
 $$
 \gamma'_t=\Ad(u^*e^{it(H\otimes1-1\otimes H^t+h)}u)
 $$                                                 
is a cocycle perturbation of $\gamma_t=\Ad\,U_t$.

Let $\phi$ be a pure ground state of $A\otimes A$ with respect to 
$\gamma'$. Then $\phi$ is a product state, i.e., 
$\phi=(\phi|B\otimes B^t)\otimes (\phi|(A\cap B')\otimes (A\cap 
B')^t)$. We consider $B$ and $B^t$ (resp. $A\cap B'\cong M_{N/K}$ 
and $(A\cap B')^t$) irreducibly acting on $\C^K$ (resp. 
$\C^{N/K}$); and then $B\otimes B^t$ on $\C^K\otimes \C^K$ (resp. 
$A\cap B'\otimes (A\cap B')^t$ on $\C^{N/K}\otimes \C^{N/K}$). 
Then there are unit vectors $\Phi_1\in \C^K\otimes \C^K$ and 
$\Phi_2\in \C^{N/K}\otimes \C^{N/K}$ such that 
 $$
 \phi(x\otimes y)=\lan x\Phi_1,\Phi_1\ran\lan y\Phi_2,\Phi_2\ran,
 \ \ x\in B\otimes B^t,\ y\in A\cap B'\otimes (A\cap B')^t.
 $$   
Let $\Phi=\Phi_1\otimes \Phi_2$ and note that 
 $$
 u^*(H\otimes 1-1\otimes H^t+h)u\Phi=
 E_0\Phi,
 $$
where $E_0$ is the minimum of the spectrum of $H\otimes 1-1\otimes 
H^t+h$.  

We may assume that $\min\Spec(H)=0$ and let $E_1=\max\Spec(H)$, 
i.e., $0,E_1\in \Spec(H)$ and $\Spec(H)\subset [0,E_1]$. From 
$H\otimes1-1\otimes H^t-\|h\|\leq H\otimes1-1\otimes H^t+h\leq 
H\otimes 1-1\otimes H^t+\|h\|$, it follows that $-E_1-\|h\|\leq 
E_0\leq -E_1+\|h\|$. Hence we have that 
 $$
 u^*(H\otimes 1-1\otimes H^t+E_1)u\Phi=(E_0+E_1)\Phi-u^*hu\Phi
 $$ 
has norm less than $2\|h\|\approx0$. Hence the distance of $u\Phi$ 
to the spectral subspace of $H\otimes1-1\otimes H^t$ corresponding 
to $[-E_1,-E_1+\eps]$ is sufficiently small for some small 
$\eps>0$ (depending on $\|h\|$ and $\|u-1\|$). Since $u\Phi\approx 
\Phi$, we thus find a unit vector $\Psi$ in the spectral subspace 
of $H\otimes1-1\otimes H^t$ corresponding to $[-E_1,-E_1+\eps]$ 
such that $\Psi\approx \Phi\equiv\Phi_1\otimes \Phi_2$. 
Specifically we may assume that 
 $$
 \lan\Phi,\Psi\ran=\Re\lan \Phi,\Psi\ran>1-\eps.
 $$  
We should note that $\Psi$ belong to the spectral subspace of 
$H\otimes 1$ corresponding to $[0,\eps]$.

Let $\calP_1$ (resp. $\calP_2$) be a maximal set of mutually 
orthogonal one-dimensional projections in $B^t$ (resp. $(A\cap 
B')^t$). Since $\sum_{p\in \calP_1}p=1$ and $\sum_{p\in 
\calP_2}p=1$, we have that 
 $$
 \sum_{p_1\in\calP_1,p_2\in\calP_2} \|(1\otimes p_1p_2)\Phi\|
 \|(1\otimes p_1p_2)\Psi\|\geq \sum\Re\lan (1\otimes 
 p_1p_2)\Phi,\Psi\ran=\lan \Phi,\Psi\ran>1-\eps
 $$
and 
 $$
 \sum_{p_1\in\calP_1,p_2\in\calP_2} \|(1\otimes p_1p_2)\Phi\|
 \|(1\otimes p_1p_2)\Psi\|\|\Phi_{p_1p_2}-\Psi_{p_1p_2}\|^2
 <2\eps,
 $$
where $\Phi_{p_1p_2}$ (resp. $\Psi_{p_1p_2}$) is the unit vector 
$c(1\otimes p_1p_2)\Phi$ (resp. $c(1\otimes p_1p_2)\Psi$) with 
normalization constant $c>0$. (If $(1\otimes p_1p_2)\Phi=0$ (resp.
$(1\otimes p_1p_2)\Psi=0$), we can disregard it.) Hence there must 
be $p_1\in \calP_1$ and $p_2\in\calP_2$ such that $(1\otimes 
p_1p_2)\Phi\not=0$, $(1\otimes p_1p_2)\Psi\not =0$, and 
 $$
 \| \Phi_{p_1p_2}-\Psi_{p_1p_2}\|^2 <3\eps
 $$                                       
(if $\eps\leq1/3$). We define a representation $\pi_0$ of $A$ on 
$\Hil_{\pi_0}=\Ran(1\otimes p_1p_2)$ by $x\mapsto x\otimes 
p_1p_2$. Since $\pi_0$ is irreducible and 
$\Phi_{p_1p_2},\Psi_{p_1p_2}\in \Hil_{\pi_0}$, there is a unitary 
$v\in A$ such that $v\approx 1$ (depending on $\eps^{1/2}$) and 
$\pi_0(v)\Phi_{p_1p_2}=\Psi_{p_1p_2}$. Furthermore, by the choice 
of $\Psi$, there is a $k\in A_{sa}$ such that $k\approx0$ 
(depending on $\eps$) and 
$\pi_0(H)\Psi_{p_1p_2}=-\pi_0(k)\Psi_{p_1p_2}$. Hence we have that 
 $$
 \pi(v^*(H+k)v)\Phi_{p_1p_2}=0.
 $$

Note that the state $\phi'$ of $A\cong A\otimes 1$ defined through 
$\pi_0$ by the unit vector $\Phi_{p_1p_2}=c(1\otimes 
p_1)\Phi_1\otimes (1\otimes p_2)\Phi_2$ with normalization 
constant $c>0$ is a pure product state with respect to $A=B\otimes 
(A\cap B')$. Since 
$\Ad(v^*e^{it(H+k)}v)(B)\stackrel{\delta'}{\subset} B$ with 
$\delta'=\delta+2\|v-1\|+2\|k\|$ for $t\in[0,1]$, we again reach 
the situation where the following lemma is applicable, which is 
already used once before.

\begin{lem}\label{A2} 
For any $\eps>0$ there exists a $\delta>0$ satisfying the 
following condition: Let $A=M_N$ for some $N\in \N$, $B$ a unital 
$*$-subalgebra of $A$ with $B\cong M_K$ for some $K$ dividing $N$, 
$H\in A_{sa}$, $\phi_1$ a pure state of $B$, and $\phi_2$ a pure 
state of $A\cap B'$ such that 
 $$
 e^{itH}Be^{-itH}\stackrel{\delta}{\subset} B
 $$
for $t\in [0,1]$ and the state $\phi_1\otimes\phi_2$ of 
$A=B\otimes(A\cap B')$ is left invariant under the action $\alpha:
t\mapsto \Ad\,e^{itH}$. Then there exist a $u\in \U(A)$ and $h\in 
A_{sa}$ such that $\|u-1\|<\eps$, $\|h\|<\eps$, and 
 $$
 u^*e^{it(H+h)}uBu^*e^{-it(H+h)}u=B,\ \ t\in \R.
 $$    
\end{lem} 
\begin{pf}
In the GNS representation $(\pi_\phi,\Hil_\phi,\Omega_\phi)$ of 
$A$ associated with $\phi=\phi_1\otimes\phi_2$ we define a 
projection $E$ on $\Hil_\phi$ by $ E=[\pi_\phi(B)\Omega_\phi]$. 
Here $[S]$ means either the projection onto the closed linear span 
of $S(\subset \Hil_\phi)$ or the closed linear span itself.                 
Then, by Lemma \ref{I} below, $\max_{t\in 
[0,1]}\|e^{itH}Ee^{-itH}-E\|$ is arbitrarily small depending on 
$\delta$. Hence, by Lemma \ref{A1} below,  there is a $u\in \U(A)$ 
and $h\in A_{sa}$ such that $u\approx1$, $h\approx0$, 
$\pi_\phi(u)\Omega_\phi=\Omega_\phi$, $\pi_\phi(h)\Omega_\phi=0$, 
and 
 $$                          
 \Ad\pi_\phi(u^*e^{it(H+h)}u)(E)=E.
 $$ 

Thus we may just as well assume that $\Ad\pi_\phi(e^{itH})(E)=E$ 
retaining the assumptions of the Lemma. 

Note that $E\in \pi_\phi(B)'$. Furthermore if $x\in B$, then 
$\pi_\phi\alpha_t(x)E=\pi_\phi(e^{itH}xe^{-itH})E=E\pi_\phi\alpha_t(x)$, 
which implies that $\pi_\phi\alpha_t(B)E=\pi_\phi(B)E$. Since 
$x\mapsto \pi(x)E$ is an isomorphism of $B$ onto $\pi_\phi(B)E$, 
one can find a map $\beta_t$ of $B$ into $B$ such that 
$\pi_\phi\alpha_t(x)E=\pi_\phi\beta_t(x)E,\ x\in B$. We can show 
that $t\mapsto\beta_t$ is a flow on $B$; e.g.,
$\pi_\phi(\beta_s\beta_t(x))E=\pi_\phi(\alpha_s(\beta_t(x)))E=\pi_\phi(e^{itH}) 
\pi_\phi(\beta_t(x))E\pi_\phi(e^{-itH})
=\pi_\phi(e^{itH})\pi_\phi\alpha_t(x)E\pi_\phi(e^{itH})=\pi_\phi\alpha_{s+t}(x)E,\ 
x\in B$.  

Let $e_1\in B$ be a minimal projection such that 
$\pi_\phi(e_1)\Omega_\phi=\Omega_\phi$. Since 
$\pi_\phi\alpha_t(e_1)\Omega_\phi=\Omega_\phi$ (because $\phi$ is 
$\alpha$-invariant), it follows that 
$\pi_\phi(\alpha_t(e_1))E\geq\pi_\phi(e_1)E$, which implies that 
$\pi_\phi(\alpha_t(e_1))E=\pi_\phi(e_1)E$ and $\beta_t(e_1)=e_1$. 
Since $\alpha_t(B)\stackrel{\delta}{\subset}B$ for $t\in [0,1]$, 
there is an $e_t\in B$ such that $\|\alpha_t(e_1)-e_t\|\leq\delta$ 
for $t\in [0,1]$. This implies that 
$\|\pi_\phi(e_1-e_t)E\|\leq\delta$, i.e., $\|e_1-e_t\|\leq\delta$. 
Hence we get that $\|\alpha_t(e_1)-e_1\|\leq 2\delta$ for $t\in 
[0,1]$. 

We will then argue that there is a projection $p\in A$ and a 
unitary $v\in A$ such that $p\approx e_1$, $\|[H,p]\|\approx 0$,
$v\approx 1$, $ve_1v^*=p$, $\pi_\phi(p)\Omega_\phi=\Omega_\phi$, 
and $\pi_\phi(v)E=E$.                  

Let $f$ be a non-negative $C^\infty$ function on $\R$ of compact 
support such that $\int f(t) dt=1$, and consider 
 $$
 \int f(t)\alpha_t(e_1)dt,
 $$
which is still close to $e_1$ for a small $\delta>0$.  By 
functional calculus construct a projection $p$ out of it, which 
satisfies that $\|[H,p]\|\approx0$ depending on $\int|f'(t)|dt$ 
(see the proof of \ref{A1} below). Since 
 $$
 \pi_\phi(\int f(t) \alpha_t(e_1)dt)E=\pi_\phi(e_1)E,
 $$
we get that $\pi_\phi(p)E=E\pi_\phi(p)=E\pi_\phi(e_1)$.  
 
Consider 
 $$
 z=pe_1+(1-p)(1-e_1),
 $$
which satisfies that $z\approx1$, $\pi_\phi(z)E=E\pi_\phi(z)=E$, 
and $ze_1=pz$, and construct a unitary $v$ by the polar 
decomposition of $z=|z|v$. Then it follows that $\pi_\phi(v)E=E$ 
and $ve_1v^*=p$.      

Let $h=-[H,p]p+p[H,p]=-(1-p)Hp-pH(1-p)$, which has small norm as 
we have assumed that $\|[H,p]\|\approx0$. Then we have that 
$\pi_\phi(h)E=0$ and $[H+h,p]=[H,p]-(1-p)Hp+pH(1-p)=0$. Hence we 
have reached the following situation:

The flow $t\mapsto \Ad(v^*e^{it(H+h)}v)$ leaves $e_1$ invariant 
and also the projection $E$ onto $\pi_\phi(B)\Omega_\phi$ 
invariant. Thus we may just as well assume that 
$\alpha_t=\Ad\,U_t$ leaves $e_1\in B$ invariant, in addition to 
the conditions already assumed. 
    
Recall that we have defined the flow $\beta$ on $B$ by 
$\pi_\phi(\beta_t(x))E=\pi_\phi(\alpha_t(x))E,\ x\in B$. Since 
$\alpha_t(B)\stackrel{\delta}{\subset}B$ for $t\in [0,1]$, we 
have, for $x\in B$ with $\|x\|\leq1$ and $t\in [0,1]$, an $x_t\in 
B$ such that $\|\alpha_t(x)-x_t\|\leq\delta$, which implies that 
$\|(\pi_\phi(\alpha_t(x))-\pi_\phi(x_t))E\|\leq \delta$. Hence we 
have that $\|(\pi_\phi(x_t)-\pi_\phi(\beta_t(x)))E\|\leq \delta$ 
or $\|x_t-\beta_t(x)\|\leq \delta$. Thus we obtain that 
 $$
 \|\alpha_t(x)-\beta_t(x)\|\leq 2\delta
 $$                                    
for $x\in B$ with $\|x\|\leq1$ and $t\in [0,1]$. 
                                                     
Since $\beta$ is a flow on $B\cong M_K$, there is a set $(e_{ij})$ 
of matrix units for $B$ such that $e_{11}=e_1$ and 
$\beta_t(e_{ij})=e^{it(p_i-p_j)}e_{ij}$ for some 
$p_1=0,p_2,\ldots,p_K\in \R$. We then define 
 $$
 v_t=\sum_{j=1}^Ke^{ip_jt}e_{j1}\alpha_t(e_{1j}).
 $$         
Then, since $\alpha_t(e_{11})=e_{11}$, $v_t$ is a unitary in $A$ 
and $t\mapsto v_t$ is an $\alpha$-cocycle such that 
$v_te_{11}=e_{11}$ and $\pi_\phi(v_t)E=E$. It follows that 
$\Ad\,v_t\alpha_t|B=\beta_t|B$, since 
 $$
 \Ad\,v_t\alpha_t(e_{ij})=e^{ip_it}e_{i1}\alpha_t(e_{1i})
 \alpha_t(e_{ij})\alpha_t(e_{j1})e_{1j}e^{-ip_jt}
 =e^{i(p_i-p_j)t}e_{ij}.
 $$ 

Since $\|(\alpha_t-\beta_t)|B\|\leq 2\delta$ for $t\in [0,1]$, we 
get that $\|(\Ad\,v_t^*-\id)|B\|=\|(\Ad\,v_t^*\beta_t-\beta_t)|B\| 
=\|(\alpha_t-\beta_t)|B\|\leq 2\delta$. Now we assert that 
$\max_{t\in [0,1]}\|v_t-1\|$ is small depending on $\delta$ (but 
without depending on the sizes of $A$ and $B$).

Since $v_t^*Bv_t=\alpha_t(B)\stackrel{\delta}{\subset}B$ for $t\in 
[0,1]$, there is a $w_t\in \U(A)$ such that $\|w_t-1\|\leq 
120\delta^{1/2}$ and $w_tv_t^*Bv_tw_t^*=B$ \cite{Ch}. Since 
$\|(\Ad(w_tv_t^*)-\id)|B\|\leq 240\delta^{1/2}+2\delta$, there 
exists a $u_t\in \U(B)$ such that $\|u_t-1\|$ is of the order 
$\|\Ad(w_tv_t^*)-\id\|^{1/2}$ and $\Ad(w_tv_t^*)|B=\Ad\,u_t^*|B$ 
(see 8.7.5 of \cite{Ped}). That is, we have that 
$z_t=u_tw_tv_t^*\in A\cap B'$, or 
 $$
 v_t=z_t^*u_tw_t.
 $$
Since $e_{11}=e_{11}v_t=e_{11}z_t\cdot e_{11}u_tw_t$ and 
$e_{11}u_tw_t\approx e_{11}$ (because $u_t\approx1\approx w_t$), 
it follows that $e_{11}z_t\approx e_{11}$. Since 
$\|z_t-1\|=\|(z_t-1)e_{11}\|$ (because $z_t\in A\cap B'$), it 
follows that $z_t\approx 1$. Thus we get that $v_t\approx 1$. 
Hence incorporating all the cocycle perturbations made we have 
reached the following situation: There is a cocycle $u$ with 
respect to $\alpha:t\mapsto \Ad\,e^{itH}$ on $A$ such that 
$\sup_{t\in [0,1]}\|u_t-1\|<\eps$ and $\Ad\,u_t\alpha_t(B)=B$. 
Then the conclusion will follow from the following lemma. 
\end{pf}        

\begin{lem}  \label{S}
For any $\eps>0$ there exists a $\delta>0$ satisfying the 
following condition: If $\alpha$ is a flow on a unital \cstar\ $A$ 
and $u$ is an $\alpha$-cocycle such that $\|u_t-1\|<\delta,\ t\in 
[0,1]$, there are $v\in \U(A)$ and a differentiable 
$\alpha$-cocycle $w$ such that $\|v-1\|<\eps$, $\|d 
w_t/dt|_{t=0}\|<\eps$ for $t\in [0,1]$, and 
$u_t=vw_t\alpha_t(v^*)$. 
\end{lem}  
\begin{pf}
This can be shown by using the 2 by 2 trick due to Connes (see 
\cite{Kis00}). Namely define a flow $\gamma$ on $M_2\otimes A$ by 
$\gamma_t(e_{11}\otimes x)=e_{11}\otimes \alpha_t(x)$, 
$\gamma_t(e_{12}\otimes x)=e_{12}\otimes \alpha_t(x)u_t^*$, etc. 
and note that $\gamma(e_{21}\otimes 1)=e_{21}\otimes u_t$. We 
approximate $e_{21}\otimes 1$ by $e_{21}\otimes v$ with 
$v\in\U(A)$ such that $d\gamma_t(e_{21}\otimes v)/dt$ has small 
norm. Let $w_t=v^*u_t\alpha_t(v)$, which is an $\alpha$-cocycle. 
Since 
 $$
 \gamma_t(e_{21}\otimes v)=e_{21}\otimes u_t\alpha_t(v)
 =e_{21}\otimes vw_t,
 $$                                                  
this concludes the proof. \end{pf} 
                                      
\begin{lem} \label{I}
For any $\eps>0$ there exists a $\delta>0$ satisfying the 
following condition: Let $A$ be a unital \cstar, $A_1$ a unital 
C$^*$-subalgebra of $A$ with $A_1\cong M_K$ for some $K\in\N$, and 
$\alpha$ a flow on $A$ such that 
 $$
 \sup_{t\in [0,1]}\dist(A_1,\alpha_t(A_1))<\delta.
 $$                                               
Let $\phi$ be an $\alpha$-invariant pure state of $A$ such that 
$\phi|A_1$ is pure and let 
 $$
 E=[\pi_\phi(A_1)\Omega_\phi].
 $$
If $U$ denotes the unitary flow on the Hilbert space $\Hil_\phi$ 
defined by 
 $$
 U_t\pi_\phi(x)\Omega_\phi=\pi_\phi(\alpha_t(x))\Omega_\phi,\ x\in 
 A, 
 $$ 
then it follows that 
 $$
 \sup_{t\in[0,1]}\|U_tEU_t^*-E\|<\eps.
 $$
\end{lem}   
\begin{pf}
Since $\phi|A_1$ is pure, there is a minimal projection $e\in A_1$ 
such that $\pi_\phi(e)\Omega_\phi=\Omega_\phi$. Note that 
$\pi_\phi(\alpha_t(e))\Omega_\phi=U_t\pi_\phi(e)\Omega_\phi=\Omega_\phi$. 
We will assert that $\|e-\alpha_t(e)\|$ is small. 

Let $t\in [0,1]$. There is an $x\in \alpha_t(A_1)$ such that 
$\|x-e\|<\delta$ and $x^*=x$. Since 
 $$
 \pi_\phi(\alpha_t(e)x\alpha_t(e))\Omega_\phi-\Omega_\phi
 =\pi_\phi(\alpha_t(e)x)\Omega_\phi-\pi_\phi(\alpha_t(e)e)\Omega_\phi
 =\pi_\phi(\alpha_t(e)(x-e))\Omega_\phi,
 $$ 
we have that 
$\|\pi_\phi(\alpha_t(e)x\alpha_t(e))\Omega_\phi-\Omega_\phi\|<\delta$. 
Since $\alpha_t(e)$ is a minimal projection in $\alpha_t(A_1)$, 
there is a $\lambda\in\R$ such that 
$\alpha_t(e)x\alpha_t(e)=\lambda\alpha_t(e)$. The above inequality 
shows that $|\lambda-1|<\delta$ or 
$\|\alpha_t(e)x\alpha_t(e)-\alpha_t(e)\|<\delta$. Since 
$\|x-e\|<\delta$, we then have that 
 $$
 \|\alpha_t(e)-\alpha_t(e)e\alpha_t(e)\|<2\delta.
 $$
Hence it follows that 
 $$
 \|\alpha_t(e)-\alpha_t(e)e\|^2=\|\alpha_t(e)(1-e)\alpha_t(e)\|<2\delta.
 $$
In the same way we get that $\|e-e\alpha_t(e)\|^2<2\delta$. Since 
$\|e-\alpha_t(e)\|\leq 
\|e-e\alpha_t(e)\|+\|e\alpha_t(e)-\alpha_t(e)\|$, we obtain that
 $$
 \|e-\alpha_t(e)\|<2\sqrt{2\delta}<3\delta^{1/2}.
 $$   
 
Let $w\in A_1$ be a partial isometry such that $w^*w=e$. We will 
assert that there is a partial isometry $y\in\alpha_t(A_1)$ such 
that $y^*y=\alpha_t(e)$ and $\|y-w\|$ is small.     

Let $z\in \alpha_t(A_1)$ such that $\|z-w\|<\delta$. Note that 
$\|z\alpha_t(e)-w\|\leq 
\|(z-w)\alpha_t(e)\|+\|w(\alpha_t(e)-e)\|<\delta+3\delta^{1/2}$. 
If $\delta$ is sufficiently small, then $\lambda=\|z\alpha_t(e)\|$ 
is close to $1$, as $|\lambda-1|<\delta+3\delta^{1/2}$. Then 
$y=\lambda^{-1}z\alpha_t(e)$ is a partial isometry in 
$\alpha_t(A_1)$ such that  $y^*y=\alpha_t(e)$ and 
$\|y-w\|<2\delta+6\delta^{1/2}$. The latter follows because $ 
\|y-w\|=\|\lambda^{-1}z\alpha_t(e)-w\|\leq \|z\alpha_t(e)-w\|+ 
 \|(\lambda^{-1}-1)z\alpha_t(e)\|<\delta+3\delta^{1/2}+|1-\lambda|$.
 
Let $\ol{\alpha}_t=\Ad\,U_t$ as a weakly continuous flow on 
$B(\Hil_\phi)$ and note that 
$\ol{\alpha}_t(E)\Hil_\phi=[\pi_\phi(\alpha_t(A_1))\Omega_\phi]$. 
We will assert that 
 $$
 \inf\Spec(E\ol{\alpha}_t(E)E)
 =\inf_{\xi}\|\ol{\alpha}_t(E)E\xi\|^2  
 $$
is close to $1$, where the spectrum is taken as an operator on 
$E\Hil_\phi$ and  the infimum is taken over all unit vectors 
$\xi\in E\Hil_\phi=[\pi_\phi(A_1)\Omega_\phi]$. Note that this 
infimum is obtained as 
 $$
 \inf_{w}\sup_{y}|\lan\pi_\phi(y)\Omega_\phi,\pi_\phi(w)\Omega_\phi\ran|^2,
 $$
where $w$ runs over all $w\in A_1$ with $w^*w=e$ and $y$ runs over 
all  $y\in \alpha_t(A_1)$ with $y^*y=\alpha_t(e)$. For each $w\in 
A_1$ with $w^*w=e$, we choose $y\in \alpha_t(A_1)$ with 
$y^*y=\alpha_t(e)$ such that $\|w-y\|<2\delta+6\delta^{1/2}$. 
Since $\|y^*w-e\|\leq \|(y^*-w^*)w\|\leq \|y-w\|$ and 
 $$
 \lan\pi_\phi(y)\Omega_\phi,\pi_\phi(w)\Omega_\phi\ran
 =1-\lan \Omega_\phi,\pi_\phi(e-y^*w)\Omega_\phi\ran,
 $$
We have that 
 $$
 |\lan\pi_\phi(y)\Omega_\phi,\pi_\phi(w)\Omega\ran|>1-(2\delta+6\delta^{1/2}).
 $$
Hence we get that
 $$
  \inf\Spec(E\ol{\alpha}_t(E)E)>1-(4\delta+12\delta^{1/2}).
 $$
 
Thus, if $t\in[0,1]$, we have that 
$\|E\ol{\alpha}_t(E)E-E\|<\delta_1$, where 
$\delta_1=4\delta+12\delta^{1/2}$. In the same way we have that 
$\|\ol{\alpha}_t(E)E\ol{\alpha}_t(E)-\ol{\alpha}_t(E)\|<\delta_1$. 
Then we get $\|E-\ol{\alpha}_t(E)\|<3\delta_1^{1/2}$ as before. 
Since $\delta_1^{1/2}\approx 2\sqrt{3}\delta^{1/4}$, this 
concludes the proof. 
\end{pf}

\begin{lem} \label{A1}
For any $\eps>0$ there exists a $\delta>0$ satisfying the 
following condition: Let $\Hil$ be a Hilbert space, $E$ a 
projection on $\Hil$, and $H$ a self-adjoint operator on $\Hil$ 
such that 
 $$
 \|e^{itH}Ee^{-itH}-E\|<\delta
 $$
for $t\in [0,1]$. Then there exist a unitary $u$ on $\Hil$ and a 
bounded self-adjoint operator $h$ on $\Hil$ such that 
$\|u-1\|<\eps$, $\|h\|<\eps$, and 
 $$
 u^*e^{it(H+h)}uEu^*e^{-it(H+h)}u=E.
 $$         
Moreover if $\Omega\in\Hil$ is a unit vector such that 
$E\Omega=\Omega$ and $e^{itH}\Omega=\Omega$, the above $u$ and $h$ 
can be chosen so that $u\Omega=\Omega$ and $h\Omega=0$. 
\end{lem}  
\begin{pf}
Let $f$ be a non-negative $C^\infty$ function on $\R$ of compact 
support such that $\int f(t)dt=1$. Then for any projection $E$ and 
a self-adjoint operator $H$ on $\Hil$, we define, for $n\in \N$, 
 $$
 E_n=\frac{1}{n}\int f(t/n)\Ad\,e^{itH}(E)dt
 $$
and estimate 
 $$
 \|[iH,E_n]\|\leq \frac {1}{n}\int |f'(t)|dt,
 $$                                          
which is close to zero depending on $n$ only. If $\delta$ is 
sufficiently small, i.e., $\|E_n-E\|\approx0$, then we can define 
a projection $F_n$ by functional calculus by 
 $$
 F_n=\frac{1}{2\pi i}\oint_C (E_n-z)^{-1}dz,
 $$
where $C$ is the path $|z-1|=1/2$ (see, e.g., \cite{Sak}). Here we 
may assume that the distance between $C$ and $\Spec(E_n)$ is 
greater than $1/4$. Since $F_n$ is close to $E_n$, this is also 
close to $E$ depending on $\|E-E_n\|$. Since 
 $$
 [iH,F_n]=\frac{-1}{2\pi i}\oint_C (E_n-z)^{-1}[iH,E_n](E_n-z)^{-1}dz,
 $$
the norm of $[iH,F_n]$ is smaller than  $16\|[iH,E_n]\|$.  Since 
we may define $u$ to be the unitary obtained from the polar 
decomposition of $z=F_nE+(1-F_n)(1-E)=1-(1-F_n)E-F_n(1-E)\approx1$ 
and $h$ as $-(1-F_n)HF_n-F_nH(1-F_n)\approx0$, this completes the 
proof of the first part. 

If $\Omega$ is a unit vector such that $E\Omega=\Omega$ and 
$H\Omega=0$, it follows that $E_n\Omega=\Omega$, 
$F_n\Omega=\Omega$, $z\Omega=\Omega$, $h\Omega=0$, i.e., the last 
conditions follow automatically.     
\end{pf}                                                             

\section{UHF algebras}  

The previous result \ref{A} can be extended to approximate UHF 
flows. 

\begin{prop} \label{M3}
For any $\eps>0$ there exists a $\delta>0$ satisfying the 
following condition: If $\alpha$ is a flow on a UHF algebra $A$ 
and $(A_n)$ is an increasing sequence $(A_n)$ of finite 
dimensional unital C$^*$-subalgebras of $A$ such that 
$\ol{\bigcup_nA_n}=A$, $A_n\cong M_{K_n}$ for all $n$, and 
$\delta_n\equiv\sup_{t\in [0,1]}\dist(\alpha_t(A_n),A_n)\ra0$ with 
$\delta_1\leq \delta$, then there is an $\alpha$-cocycle $u$ such 
that $\sup_{t\in [0,1]} \|u_t-1\|<\eps$ and 
$\Ad\,u_t\alpha_t(A_n)=A_n$ for $n=1$ and infinitely many $n$. 
\end{prop} 
\begin{pf}
Let $\tau$ denote the tracial state of the UHF algebra $A$. We 
denote by $(\pi_\tau,\Hil_\tau,\Omega_\tau)$ the GNS 
representation associated with $\tau$ and recall that the 
canonical conjugation operator $J$ is defined by 
 $$
 J\pi_\tau(x)\Omega_\tau=\pi_\tau(x^*)\Omega_\tau, \ \ x\in A.
 $$
By using the fact that $J\pi_\tau(A)''J=\pi_\tau(A)'$, we define 
an irreducible representation $\rho$ of $A\otimes A^{op}$ in 
$\Hil_\tau$ by 
 $$
 \rho(x\otimes y)=\pi_\tau(x)J\pi_\tau(y^*)J,\ \ x\otimes y\in 
 A\otimes A^{op},
 $$                                          
where $A^{op}$ is the opposite \cstar\ of $A$; i.e., $A^{op}=A$ as 
a Banach space with the same involution and the new product 
$x\circ y=yx$. We check that $\rho$ is indeed a representation of 
$A\otimes A^{op}$ by  $\rho(x\otimes y)\rho(a\otimes 
b)=\pi_\tau(xa)J\pi_\tau(y^*b^*)J=\pi_\tau(xa)J\pi_\tau((y\circ 
b)^*)J=\rho(xa\otimes y\circ b)$ and note that $A\otimes A^{op}$ 
is a UHF algebra with dense $\bigcup_nA_n\otimes A_n$. We define a 
state $\omega$ of $A\otimes A^{op}$ by 
$\omega(z)=\lan\rho(z)\Omega_\tau,\Omega_\tau\ran$. Since 
$\rho(A_n\otimes A_n)\Omega_\tau=\pi_\tau(A_n)\Omega_\tau$ is 
$K_n^2$-dimensional, we know that $\omega|A_n\otimes A_n$ is pure 
for all $n$. (If $(e_{ij})$ is a family of matrix units for $A_n$, 
then $p=K_n^{-1}\sum_{ij}e_{ij}\otimes e_{ji}$ is a minimal 
projection in $A_n\otimes A_n$ such that 
$\rho(p)\Omega_\tau=\Omega_\tau$.) 

We define a unitary flow $U$ in $\Hil_\tau$ by 
 $$
 U_t\pi_\tau(x)\Omega_\tau=\pi_\tau\alpha_t(x)\Omega_\tau,\ \ x\in 
 A.
 $$
Since $U_t\rho(x\otimes 
y)\Omega_\tau=\pi_\tau\alpha_t(x)\pi_\tau\alpha_t(y)\Omega_\tau= 
\rho(\alpha_t(x)\otimes \alpha_t(y))\Omega_\tau$, we get that 
$\Ad\,U_t\rho=\rho\circ(\alpha_t\otimes \alpha_t)$. Note that 
$\omega\circ(\alpha_t\otimes \alpha_t)=\omega$. 

If $\delta_n$ is sufficiently small and $t\in[0,1]$, there is a 
unitary $v_n\in A$, by \cite{Ch}, such that $\|v_n-1\|\leq 
120\delta_n^{1/2}$ and $v_n\alpha_t(A_n)v_n^*=A_n$. Hence we get 
that 
 $$
 \Ad(v_n\otimes v_n^*)(\alpha_t\otimes \alpha_t)(A_n\otimes A_n)    
 =A_n\otimes A_n.
 $$        
This implies that $\dist((\alpha_t\otimes \alpha_t)(A_n\otimes 
A_n),A_n\otimes A_n)\leq 480\delta_n^{1/2}$ for $t\in[0,1]$. 

Hence we shall first show the following weaker version of this 
proposition. 
\end{pf}
                    
\begin{lem}   \label{B}
For any $\eps>0$ there exists a $\delta>0$ satisfying the 
following condition: If $\alpha$ is a flow on a UHF algebra $A$, 
$\phi$ is an $\alpha$-invariant pure state of $A$, and $(A_n)$ is 
an increasing sequence $(A_n)$ of unital finite dimensional 
C$^*$-subalgebras of $A$ such that $\ol{\bigcup_nA_n}=A$, 
$A_n\cong M_{K_n}$ for all $n$, $\phi|A_n$ is pure for all $n$, 
and $\delta_n\equiv\sup_{t\in [0,1]}\dist(\alpha_t(A_n),A_n)\ra0$ 
as $n\ra\infty$ with $\delta_1\leq \delta$, then there is an 
$\alpha$-cocycle $u$ such that $\|u_t-1\|<\eps$ for all 
$t\in[0,1]$ and $\Ad\,u_t\alpha_t(A_n)=A_n$ for $n=1$ and 
infinitely many $n$. 
\end{lem}      

We will prove the above lemma by induction. What we need is the 
following lemma for a unital \cstar\ $A$. 

\begin{lem}  \label{C}
For any $\eps>0$ there is a $\delta>0$ satisfying the following 
condition. If $\alpha$ is a flow on a unital \cstar\ $A$, $\phi$ 
is an $\alpha$-invariant pure state of $A$, and $A_1$ is a unital 
C$^*$-subalgebra of $A$ such that $A_1\cong M_{K}$ for some $K$, 
$\phi|A_1$ is pure, and 
 $$\sup_{t\in 
[0,1]}\dist(\alpha_t(A_1),A_1)<\delta, 
 $$ 
then there is an $\alpha$-cocycle $u$ such that $\|u_t-1\|<\eps$ 
for all $t\in[0,1]$, $\phi\Ad\,u_t\alpha_t=\phi$, and 
$\Ad\,u_t\alpha_t(A_1)=A_1$. 
\end{lem}    

\begin{pf}
In the GNS representation $(\pi_\phi,\Hil_\phi,\Omega_\phi)$ 
associated with $\phi$, we define a unitary flow $U$ by 
$U_t\pi_\phi(x)\Omega_\phi=\pi_\phi(\alpha_t(x))\Omega_\phi,\ x\in 
A$ and denote by $H$ the self-adjoint generator of $U$: 
$U_t=e^{itH}$. 

We define a projection $E_1$ in $\Hil_\phi$ by 
 $$
 E_1=[\pi_\phi(A_1)\Omega_\phi].
 $$                      
Note that $E_1$ is a projection in $\pi_\phi(A_1)'$ of rank $K$ 
and $A_1\ni x\mapsto \pi_\phi(x)E_1$ is an irreducible 
representation of $A_1$. As in the proof of \ref{A}, 
$\|U_tE_1U_t^*-E_1\|$ is close to zero depending on $\delta$. Then 
we find a projection $F$ of finite rank such that $E_1\approx F$, 
$F\Omega_\phi=\Omega_\phi$, $\|[H,F]\|\approx0$. We also find a 
unitary $Z$ on the subspace $L$ spanned by $E_1\Hil_\phi$ and 
$F\Hil_\phi$ such that $Z\approx1$, $Z\Omega_\phi=\Omega_\phi$, 
and $F=ZE_1Z^*$. By using Kadison's transitivity \cite{Kad,Ped}, 
we find a $u\in\U(A)$ such that $\pi_\phi(u)|L=Z|L$ and 
$\|u-1\|<2\|Z-1\|$.  We also find an $h\in A_{sa}$ such that 
$\pi_\phi(h)=-(1-F)HF-FH(1-F)$ on the subspace spanned by 
$F\Hil_\phi$ and $(1-F)HF\Hil_\phi$, and $\|h\|<2\|(1-F)HF\|$. 
Then it follows that $\pi_\phi(u)\Omega_\phi=\Omega_\phi$,  
$\pi_\phi(h)\Omega_\phi=0$, $\pi_\phi(u)E_1\pi_\phi(u^*)=F$,  and 
$[H,F]=-[\pi_\phi(h),F]$ as well as $u\approx1$ and $h\approx0$. 
Note that all the estimates depend on $\delta$ but not on the size 
of $A_1$ and that 
$\Ad(\pi_\phi(u)^*e^{it(H+\pi_\phi(h))}\pi_\phi(u))E_1=0$. If 
$e_1$ denotes the minimal projection of $A_1$ such that 
$\pi_\phi(e_1)\Omega_\phi=\Omega_\phi$, we may also suppose that 
$(\delta_\alpha+\ad(ih))\Ad\,u(e_1)=0$ (see the proof of \ref{A}). 

Let $\alpha'_t=\Ad\,u\,e^{it(\delta_\alpha+\ad\,ih)}\Ad\,u^*$, 
which is a small cocycle perturbation of $\alpha$. We then define 
a flow $\beta$ on $A_1$ by                      
 $$
 \pi_\phi(\alpha'_t(x))E_1=\pi_\phi(\beta_t(x))E_1,\ \ x\in A_1.
 $$
It follows from the commutativity of 
$\pi_\phi(u)^*e^{it(H+\pi_\phi(h))}\pi_\phi(u)$ and $E_1$ that 
$\beta$ is indeed a flow. Then we derive that 
$\|(\alpha'_t-\beta_t)|A_n\|\leq 2\delta'$, where 
$\delta'=\sup_{t\in [0,1]}\dist(\alpha_t'(A_1),A_1)$ which is 
small depending on the original $\delta\geq 
\sup_{t\in[0,1]}\dist(\alpha_t(A_1),A_1)$. We choose a family 
$(e_{ij})$ of matrix units for $A_1$ such that 
$\beta_t(e_{i1})=e^{ip_it}e_{i1}$ for some 
$p_1=0,p_2,\ldots,p_{K}$ in $\R$ with $e_{11}=e_n$. Noting that 
$e_{11}=\alpha_t'(e_{11})$, we define an $\alpha'$-cocycle 
 $$
 v_t=\sum_ie^{ip_it}e_{i1}\alpha'_t(e_{1i}).
 $$
We show that $\pi_\phi(v_t)\Omega_\phi=\Omega_\phi$, $v_t\approx1$ 
(depending on $\delta'$; see the proof of \ref{A2}), and 
$\Ad\,v_t\alpha'_t|A_1=\beta_t|A_1$. Thus 
$\alpha_t''=\Ad\,v_t\alpha_t'$, which is a small cocycle 
perturbation of $\alpha$, satisfies the required condition. Note 
that the estimate of how far $\alpha''$ is from $\alpha$ does not 
depend on the size of $A_1$, thanks to the estimate in 5.3 of 
\cite{Ch}. 
\end{pf} 

To prove Lemma \ref{B} we apply the above lemma inductively. By 
choosing $\delta$ small enough, we guarantee the above lemma 
applies to $A_1$; thus we find an $\alpha$-cocycle $u^1$ such that 
$\max_{t\in[0,1]}\|u^1_t-1\|<\eps/2$, 
$\alpha^1_t=\Ad\,u^1_t\alpha_t$ fixes $A_1$, and 
$\phi\alpha^1_t=\phi$. To proceed to the next step, we just note 
that 
 $$
 \delta_n\equiv\sup_{t\in[0,1]}(\alpha^1_t(A_n\cap A_1'),A_n\cap A_1')\ra0
 $$
as $n\ra\infty$. We find $n_2>1$ such that $\delta_{n_2}$ is 
sufficiently small so that we find an $\alpha^1$-cocycle $u^2$ in 
$A\cap A_1'$ such that $\max_{t\in [0,1]}\|u^2_t-1\|<\eps/4$, 
$\alpha^2_t=\Ad\,u^2_t\alpha^1_t$ fixes $A_{n_2}$, and 
$\phi\alpha^2_t=\phi$. We repeat this process inductively. 

\begin{lem} \label{D}                                                   
For any $\eps>0$ there exists a $\delta>0$ satisfying the 
following condition:  Let $\alpha$ be a flow on a unital \cstar\ 
$A$ and $A_1$ a unital C$^*$-subalgebra of $A$ such that 
$A_1\equiv M_K$ for some $K$. Let $u$ be an 
$\alpha\otimes\alpha$-cocycle, where $t\mapsto 
\alpha_t\otimes\alpha_t$ is a flow on $A\otimes A^{op}$, such that 
$\max_{t\in [0,1]}\|u_t-1\|<\delta$ and $\Ad\,u_t(\alpha_t\otimes 
\alpha_t)$ fixes $A_1\otimes A_1$. Let $\phi$ be a pure ground 
state of $A\otimes A^{op}$ with respect to the flow $t\mapsto 
Ad\,u_t(\alpha_t\otimes\alpha_t)$ (such that $\phi|A_1\otimes A_1$ 
is pure). Then there is an $\alpha$-cocycle $v$ such that 
$\max_{t\in[0,1]}\|v_t-1\|<\eps$ and $\Ad\,v_t\alpha_t$ fixes 
$A_1$. 
\end{lem}    
\begin{pf}
We may suppose that $u$ is given as 
$wv_t(\alpha_t\otimes\alpha_t)(w^*)$, where $w\in \U(A\otimes 
A^{op})$ and $v$ is a differentiable $\alpha\otimes\alpha$-cocycle 
with $ih=dv_t/dt|_{t=0}$ such that $\|w-1\|<\delta$ and 
$\|h\|<\delta$ (see \ref{S}). 

Let $(\pi_\phi,\Hil_\phi,\Omega_\phi)$ be the GNS representation 
of $A\otimes A^{op}$ associated with $\phi$. Since 
$\phi_1=\phi|A_1\otimes A_1$ and $\phi_2=\phi|(A\cap 
A_1')\otimes(A\cap A_1')$ are pure, 
$(\pi_\phi,\Hil_\phi,\Omega_\phi)$ is identified with 
$(\pi_{\phi_1}\otimes \pi_{\phi_2}, 
\Hil_{\phi_1}\otimes\Hil_{\phi_2}, \Omega_{\phi_1}\otimes 
\Omega_{\phi_2})$, where 
$(\pi_{\phi_i},\Hil_{\phi_i},\Omega_{\phi_i})$ is the GNS triple 
for $\phi_i$. 

We define a unitary flow $U$ in $\Hil_\phi$ by 
 $$
 U_t\pi_\phi(x)\Omega_\phi=\pi_\phi\Ad\,u_t(\alpha_t\otimes\alpha_t)(x)\Omega_\phi,
 \ \ x\in A\otimes A^{op}.
 $$             
Let $H$ be the generator of $U$. Since $\phi$ is a ground state, 
we have that $H\geq0$ and $H\Omega_\phi=0$. Let 
$H_0=\pi_\phi(w^*)H\pi_\phi(w)-\pi_\phi(h)$, which is a 
self-adjoint operator in $\Hil_\phi$ with the domain 
$\D(H_0)=\pi_\phi(w^*)\D(H)$ and let $E_0$ be the spectral measure 
of $H_0$. Then, since $\Ad\,u_t\circ(\alpha_t\otimes\alpha_t) 
=\Ad\,w\circ\Ad\,v_t\circ(\alpha_t\otimes\alpha_t)\circ\Ad\,w^*$, 
we have that 
 $$
 \Ad\,e^{itH_0}\pi_\phi(x\otimes y)=\pi_\phi(\alpha_t(x)\otimes \alpha_t(y)),
 \ \ x,y\in A.
 $$                                                                                 
Since $H_0\geq -\|h\|>-\delta$ and 
$\|H_0\pi_\phi(w^*)\Omega_\phi\|=\|\pi_\phi(hw^*)\Omega_\phi\|<\delta$, 
we should note that 
$E_0[-\delta,\delta^{1/2}]\pi_\phi(w^*)\Omega_\phi$ has norm close 
to $1$, or more concretely, 
$\|[E_0(\delta^{1/2},\infty)\pi_\phi(w^*)\Omega_\phi\|<\delta^{1/2}$. 

Let $\pi$ denote the representation of $A$ in $\Hil_\phi$ defined 
by $\pi(x)=\pi_\phi(x\otimes 1), \ x\in A$. Since $\pi_\phi$ is 
irreducible, $\pi(A)''$ is a factor. Since $H_0$ is bounded below, 
the Borchers' theorem \cite{Sak} tells us that there is a unitary 
flow $V$ in $\pi(A)''$ such that 
$V_t\pi(x)V_t^*=U_t\pi(x)U_t^*=\pi\alpha_t(x),\ x\in A$. Let $H_1$ 
be the generator of $V$ and let $E_1$ be the spectral measure of 
$H_1$. Note that $H_1$ is bounded below and is unique up to 
constant. 

Let $P'=[\pi(A)E_0[-\delta,\delta^{1/2}]\Hil_\phi]$ and 
$P=[\pi(A)'E_0[-\delta,\delta^{1/2}]\Hil_\phi]$. Note that $P'\in 
\pi(A)'\cap U'$, $P\in \pi(A)''\cap U'$, and $PP'\geq 
E_0[-\delta,\delta^{1/2}]$. Then we have: 

\begin{lem} $PP'\leq E_0[-\delta,\delta+2\delta^{1/2}]$.
\end{lem}
\begin{pf} 
This is proved in \cite{Kis}, but we give a proof here. 

Let $\lambda\in \Spec(UPP')$. Since 
$[PP'\Hil_\phi]=[P\pi(A)PE_0[-\delta, \delta^{1/2}]\Hil_\phi]$, 
for any $\eps'>0$ there is a $Q\in P\pi(A)P$ and $\xi\in 
E_0[-\delta,\delta^{1/2}]\Hil_\phi$ such that $Q\xi\not=0$ and 
$\Spec_U(Q\xi)\subset (\lambda-\eps',\lambda+\eps')$, where 
$\Spec_U(\zeta)$ is the $U$-spectrum of $\zeta$, meaning the least 
closed subset $B$ of $\R$ with $E_0(B)\zeta=\zeta$, for 
$\zeta\in\Hil_\phi$ (see \cite{Ped}). We may assume that 
$\Spec_{\ol{\alpha}}(Q)\subset 
(\lambda-\eps'-\delta^{1/2},\lambda+\eps'+\delta)$, where 
$\ol{\alpha}$ is the (weakly continuous) flow $t\mapsto \Ad\,U_t$ 
on $B(\Hil_\phi)$ and $\Spec_{\ol{\alpha}}(Q)$ is the 
$\ol{\alpha}$-spectrum of $Q\in B(\Hil_\phi)$. Since 
$Q^*E_0[-\delta,\delta^{1/2}]\not=0$, there is an $\eta\in 
E_0[-\delta,\delta^{1/2}]\Hil_\phi$ such that $Q^*\eta\not=0$. 
Since $\Spec_U(Q^*\eta)\subset 
\Spec_{\ol{\alpha}}(Q^*)+\Spec_U(\eta)\subset 
(-\lambda-\eps'-2\delta,-\lambda+\eps'+2\delta^{1/2})\cap 
\Spec(UPP')\not=\emptyset$ or 
 $$
 \lambda\in -\Spec(UPP')+(-\eps'-2\delta,\eps'+2\delta^{1/2}).
 $$
Since $\eps'>0$ is arbitrary, we get that
 $$
 \lambda\in -\Spec(UPP')+[-2\delta,2\delta^{1/2}]\subset (-\infty,\delta+2\delta^{1/2}].
 $$
Since $\lambda\geq -\delta$, this concludes the proof. 
\end{pf}

Hence it follows that for a small $\delta>0$,
 $$
 \Spec_{\ol{\alpha}}(P\pi(A)''PP')\subset [-2\delta-2\delta^{1/2},2\delta+2\delta^{1/2}]
 \subset [-3\delta^{1/2},3\delta^{1/2}].
 $$              
Since $P\pi(A)''P\ni Q\mapsto QP'\in P\pi(A)''PP'$ is an 
isomorphism intertwining $\ol{\alpha}$, it also follows that 
$\Spec_{\ol{\alpha}}(P\pi(A)''P)\subset 
[-3\delta^{1/2},3\delta^{1/2}]$ and that 
$\max\Spec(H_1P)-\min\Spec(H_1P)\leq 3\delta^{1/2}$. Hence 
adjusting a constant to $H_1$ we may suppose that  
 $$
 P\leq E_1[0,3\delta^{1/2}].
 $$                                                
Then it automatically follows that $H_1\geq 
-\delta^{1/2}-\delta>-2\delta^{1/2}$. Because if $E_1(-\infty, 
-\delta^{1/2}-\delta')\not=0$ for some $\delta'>\delta$, there is 
a non-zero $Q\in \pi(A)''$ such that 
$Q=E_1(-\infty,-\delta^{1/2}-\delta')QP$ as $\pi(A)''$ is a 
factor. Then it follows that $\Spec_{\ol{\alpha}}(Q)\subset 
(-\infty, -\delta^{1/2}-\delta']$. Since there is a $\xi\in 
E_0[-\delta,\delta^{1/2}]\Hil_\phi$ such that $Q\xi\not=0$, we 
reach the contradiction that $\Spec_U(Q\xi)\subset 
(-\infty,-\delta']\subset (-\infty,-\delta)$. 

Since $P\geq E_0[-\delta,\delta^{1/2}]$, we have that 
$E_0[-\delta,\delta^{1/2}]\leq E_1[0,3\delta^{1/2}]$. 
               
Let $\pi_1$ (resp. $\pi_2$) denote the representation of $A_1$ 
defined by $\pi_1(x)=\pi_{\phi_1}(x\otimes1),\ x\in A_1$ (resp. 
$\pi_2(x)=\pi_{\phi_2}(x\otimes 1),\ x\in A\cap A_1'$). Note that 
$\pi=\pi_1\otimes \pi_2$ while $A=A_1\otimes (A\cap A_1')$. Let 
$F=[\pi(A)\Omega_\phi]=[\pi_1(A_1)\Omega_{\phi_1}]\otimes 
[\pi_2(A\cap A_1' )\Omega_{\phi_2}]\in \pi(A)'$. We will restrict 
the representation $\pi$ to the cyclic subspace $F\Hil_\phi$ 
below. 

Since $\|(1-F)\pi_\phi(w^*)\Omega_\phi\|<\delta$ and 
 $$
 \|(1-E_1[0,3\delta^{1/2}])\pi_\phi(w^*)\Omega_\phi\|\leq 
 \|E_0(\delta^{1/2},\infty)\pi_\phi(w^*)\Omega_\phi\| 
 <\delta^{1/2},
 $$
we have that 
$\|FE_1[0,3\delta^{1/2}]\pi_\phi(w^*)\Omega_\phi-\Omega_\phi\| 
<2\delta+\delta^{1/2}$. Let 
$\Psi=cFE_1[0,3\delta^{1/2}]\pi_\phi(w^*)\Omega_\phi$ have norm 
one with $c>0$. Then we may suppose that 
$\|\Psi-\Omega_\phi\|<2\delta^{1/2}$ (for a small $\delta>0$).       

Note that $F$ is given as $F_1\otimes F_2$, where 
$F_1=[\pi_1(A_1)\Omega_{\phi_1}]$ and $F_2=[\pi_2(A\cap 
A_1')\Omega_{\phi_2}]$. We choose a maximal abelian 
W$^*$-subalgebra $C_1$ (resp. $C_2$) of $F_1\pi_1( A_1)'F_1$ 
(resp. $F_2\pi_2(A\cap A_1')'F_2$) so that $C=C_1\otimes C_2$ is 
maximal abelian in $F\pi(A)'F$, where $C_1$ is finite-dimensional. 
We apply the decomposition theory to the cyclic representation 
$F\pi$ with the cyclic vector $\Omega_\phi=\Omega_{\phi_1}\otimes 
\Omega_{\phi_2}$ with respect to this maximal abelian 
W$^*$-subalgebra $C$. Let $K_i$ be the character space of $C_i$ 
and let $\nu_i$ be the probability measure on $K_i$ defined by 
$Q\in C_i\mapsto \lan Q\Omega_{\phi_i},\Omega_{\phi_i}\ran$. Then 
one can express $(F\pi,F\Hil_\phi,\Omega_\phi)$ as a direct 
integral over $(K_1\times K_2,\nu_1\otimes\nu_2)$; e.g., 
 $$
 F\Hil_\phi=\int_{K_1\times K_2}^{\oplus} 
 \Hil_1(s_1)\otimes \Hil_2(s_2)d\nu_1(s_1)d\nu_2(s_2),
 $$                                            
where $(\pi_i(s),\Hil_i(s),\xi_i(s))$ is the representation 
associated with $s\in K_i$. Note that $\pi_i(s)$ is irreducible 
for almost all $s\in K_i$. 

Recall that $\Psi\in F\Hil_\phi$ and let $\mu$ be the probability 
measure on $K_1\times K_2$ defined by $Q\in C\mapsto \lan 
Q\Psi,\Psi\ran$. Then $\mu$ is absolutely continuous with respect 
to $\nu=\nu_1\otimes\nu_2$; so we set $f=d\mu/d\nu$. Then we have 
 $$
 \Psi=\int^{\oplus}_{K_1\times K_2}\eta(s_1,s_2)f(s_1,s_2)d\nu_1(s_1)d\nu_2(s_2).
 $$                    
where $\eta(s_1,s_2)$ is a unit vector in $\Hil_1(s_1)\otimes    
\Hil_2(s_2)$. 

Since $(H_1+1)^{-1}\in \pi(A)''$, we get a measurable function 
$(s_1,s_2)\mapsto  B(s_1,s_2)\in (\pi_1\otimes \pi_2)(A)''$ such 
that 
 $$
 (H_1+1)^{-1}F=\int^{\oplus}_{K_1\times 
 K_2}B(s_1,s_2)d\nu_1(s_1)d\nu_2(s_2).
 $$
For almost all $(s_1,s_2)\in K_1\times K_2$, $B(s_1,s_2)$ is given 
as $(H_1(s_1,s_2)+1)^{-1}$, where $H_1(s_1,s_2)$ is self-adjoint 
with $H_1(s_1,s_2)\geq -2\delta^{-1/2}$. Setting 
$V_t(s_1,s_2)=e^{itH_1(s_1,s_2)}$, we get that   
 $$
 V_tF=\int^{\oplus}_{K_1\times 
 K_2}V_t(s_1,s_2)d\nu_1(s_1)d\nu_2(s_2),\ \ t\in\R.
 $$  
Since $E_1[0,3\delta^{1/2}]\Psi=\Psi$, it also follows for almost 
all $(s_1,s_2)$ that 
 $$
 E_1(s_1,s_2)[0,3\delta^{1/2}]\eta(s_1,s_2)=\eta(s_1,s_2),
 $$ 
where $E_1(s_1,s_2)$ is the spectral measure of $H_1(s_1,s_2)$, 
and that 
 $$
 \Ad\,V_t(s_1,s_2)(\pi_1(s_1)\otimes\pi_2(s_2))=
 (\pi_1(s_1)\otimes\pi_2(s_2))\alpha_t.
 $$  
                                            
Note that 
 \BE
 \|\Omega_\phi-\Psi\|
 &=&\int_{K_1\times K_2}\|\xi_1(s_1)\otimes 
 \xi_2(s_2)-\sqrt{f(s_1,s_2)}\eta(s_1,s_2)\|^2d\nu_1(s_1)d\nu_2(s_2)\\
 &=& \int (1+f(s_1,s_2)-2\sqrt{f(s_1,s_2)}\Re\lan 
 \xi_1(s_1)\otimes\xi_2(s_2),\eta(s_1,s_2))\ran d\nu_1(s_1)d\nu_2(s_2).
 \EE
Since $\|\Omega_\phi-\Psi\|^2<4\delta$, we have $(s_1,s_2)\in 
K_1\times K_2$ such that all the previous conditions are satisfied 
for this $(s_1,s_2)$ and $f(s_1,s_2)\leq1$ and 
 $$
 1+f(s_1,s_2)-2\sqrt{f(s_1,s_2)}\Re\lan 
 \xi_1(s_1)\otimes\xi_2(s_2),\eta(s_1,s_2)\ran <4\delta.
 $$
The latter two conditions imply that  
 $$
 1-\Re\lan 
 \xi_1(s_1)\otimes\xi_2(s_2),\eta(s_1,s_2)\ran
 \leq 1-\sqrt{f(s_1,s_2)}\Re\lan \xi_1(s_1)\otimes \xi_2(s_2),\eta(s_1,s_2)\ran
 <2\delta^{1/2},
 $$                                                      
which in turn implies that 
 $$
 \|\xi_1(s_1)\otimes \xi_2(s_2)-\eta(s_1,s_2)\|<4\delta^{1/2}.
 $$

Thus we find $s_1\in K_1$ and $s_2\in K_2$ such that the above 
norm estimate holds, $\pi_i(s_i)$ is irreducible, 
$\Ad\,V_t(s_1,s_2)$ implements $\alpha_t$ in 
$\pi_1(s_1)\otimes\pi_2(s_2)$, and 
 $$
 E_1(s_1,s_2)[0,3\delta^{1/2}]\eta(s_1,s_2)=\eta(s_1,s_2).
 $$
We can then find a $z\in \U(A)$ and $b\in A_{sa}$ such that 
$z\approx1$, $\pi_1(s_1)\otimes\pi_2(s_2)(z)\xi_1(s_1)\otimes 
\xi_2(s_2)=\eta(s_1,s_2)$,  $\|b\|\approx0$, and 
$(H_1(s_1,s_2)+(\pi_1(s_1)\otimes \pi_2(s_2))(b))\eta(s_1,s_2)=0$. 
Thus we have a small $\alpha$-cocycle $w$ such that the perturbed 
$\Ad\,w_t\alpha_t$ has a pure invariant state $\psi$ which is pure 
on $A_1$, i.e., $w_t=zv_t\alpha(z^*)$, where $v$ is the 
$\alpha$-cocycle defined by $dv_t/dt|_{t=0}=ib$ and $\psi(x)=\lan 
\pi_1(s_1)\otimes\pi_2(s_2)(x)\xi_1(s_1)\otimes 
\xi_2(s_2),\xi_1(s_1)\otimes\xi_2(s_2)\ran,\ x\in A$. Now we are 
in the situation where we can invoke Lemma \ref{C}. 
\end{pf}

\section{AF algebras}  

To prove Proposition \ref{Inv} we first show it in a special case 
as follows:  

\begin{lem} \label{I1}
For any $\eps>0$ there exists a $\delta>0$ satisfying the 
following condition: Let $A$ be a unital \cstar\ and let $\alpha$ 
be a flow on $A$. If $D$ is a unital finite-dimensional abelian 
C$^*$-subalgebra of $A$ such that 
 $$
 \sup_{|t|\leq 1} \|(\alpha_t-\id)|D\|<\delta,
 $$
then there is an $\alpha$-cocycle $u$ such that $\max_{|t|\leq 
1}\|u_t-1\|<\eps$ and $\Ad\,u_t\alpha_t(x)=x,\ x\in D$. 
\end{lem}
\begin{pf}
The main point here is that we can choose $\delta>0$ independently 
of the dimensionality of $D$; otherwise the lemma would be 
obvious. 

Suppose that we are given a unital \cstar\ $A$ and a flow $\alpha$ 
on $A$. 

Let $f$ be a non-negative $C^\infty$-function on $\R$ of compact 
support such that $\int f(t) dt=1$. We define a unital completely 
positive (or CP for short) map $\alpha_f$ on $A$ by 
 $$
 \alpha_f(x)=\int f(t)\alpha_t(x) dt.
 $$
Since $\delta_\alpha\alpha_f(x)=-\int f'(t)\alpha_t(x)dt$, we get 
that 
 $$
 \|\delta_\alpha\alpha_f\|\leq \int |f'(t)|dt.
 $$

Given a small $\eps>0$ we choose such an $f$ satisfying 
 $$
 \int|f'(t)|dt<\eps,
 $$
which entails $\|\delta_\alpha\alpha_f\|<\eps$. Then we choose 
$\delta>0$ such that 
 $$\delta\int f(t)(1+|t|)dt <\eps^2/2.
 $$ 

Suppose that we are given a unital abelian finite-dimensional 
C$^*$-subalgebra $D$ of $A$ such that $\sup_{|t|\leq 
1}\|(\alpha_t-\id)|D\|<\delta$, which entails that 
 $$
 \|(\alpha_f-\id)|D\|<\eps^2/2.
 $$
Then we get that 
 $$
 \alpha_f(u)\alpha_f(u)^* \geq 1-\eps^2
 $$
for any unitary $u\in D$. 

What we do in the following is to find a homomorphism $\Phi$ of 
$D$ into $A$ in a close neighborhood of $\alpha_f$ such that the 
norm  $\|\delta_\alpha\Phi\|$ is small depending on 
$\|\delta_\alpha\alpha_f\|$. To find such a $\Phi$ we will follow 
the strategy taken by Christensen for the proof of 3.3 of 
\cite{Ch1}.  

Let $(\pi,U)$ be a covariant faithful representation of 
$(A,\alpha)$. By Stinespring's theorem we get a representation 
$\rho$ of $A$ such that $\Hil_\rho\supset \Hil_\pi$, 
$[\rho(A)\Hil_\pi]=\Hil_\rho$,  and 
 $$
 \pi(\alpha_f(x))=P\rho(x)|\Hil_\pi,\ \ x\in A,
 $$
where $P$ is the projection onto $\Hil_\pi$, i.e., $P=\pi(1)$. 
This representation $\rho$ can be obtained as follows: We define 
an inner product on the algebraic tensor product $A\odot \Hil_\pi$ 
by 
 $$
 \lan \sum_i x_i\otimes \xi_i,\sum_j y_j\otimes \eta_j\ran 
 =\sum_i\sum_j \lan \alpha_f(y_j^*x_j)\xi_i,\eta_j\ran.
 $$
We then define a representation $\rho$ of $A$ on $A\odot 
\Hil_\phi$ by 
 $$
 \rho(a)\sum_i x_i\otimes \xi_i=\sum_i ax_i\otimes \xi_i.
 $$
Then $\Hil_\rho$ is defined as the completion of the quotient of 
$A\odot \Hil_\pi$ by the linear subspace consisting of 
$\sum_ix_i\otimes \xi_i$ with $\|\sum_ix_i\otimes \xi_i\|=0$ and 
$\rho$ naturally induces a representation of $A$ on $\Hil_\rho$, 
which we will also denote by $\rho$.  We regard $\Hil_\pi$ as a 
subspace of $\Hil_\rho$ by mapping $\Hil_\pi$ into $A\odot 
\Hil_\pi$ by $ \xi\mapsto 1\otimes \xi$, which we can easily check 
is isometric. For each $t\in \R$ we define an operator $V_t$ on 
$A\odot \Hil_\pi$ by 
 $$
 V_t \sum_i x_i\otimes \xi_i=\sum_i\alpha_t(x_i)\otimes U_t\xi_i.
 $$
Then $V=(V_t,\ t\in\R)$ induces a unitary flow on $\Hil_\rho$, 
which is denoted by the same symbol. Then we check that 
$V_tP=PV_t$, $U_t=V_tP$, and $V_t\rho(x)V_t^*=\rho\alpha_t(x),\ 
x\in A$. 

Let $B$ be the \cstar\ generated by $\rho(A)$ and $P$; then it 
follows that $PBP\subset \pi(A)$. We define a flow $\beta$ on $B$ 
by $ \beta_t=\Ad\,V_t|B$, which satisfies that 
$\beta_t\rho(x)=\rho\alpha_t(x),\ x\in A$, $\beta_t(P)=P$, and 
$\beta_t\pi(x)=\pi\alpha_t(x)$ for $x\in A$ with $\pi(x)\in PBP$. 

By mimicking the proof of 3.3 of Christensen \cite{Ch1}, we define 
 $$
 R=\int_{\U(D)}\rho(u)P\rho(u)^* d\nu(u)\in \rho(D)'\cap B,
 $$
where $\nu$ is normalized Haar measure on the compact group 
$\U(D)$. Since 
$\|\rho(u)P-P\rho(u)\|=\|(1-P)\rho(u)P-P\rho(u)(1-P)\|$ is equal 
to
 $$
 \max_{z=u,u^*} \|P\rho(z)(1-P)\rho(z)^*P\|^{1/2}
 =\max_{z=u,u^*}\|1-\alpha_f(z)\alpha_f(z)^*\|^{1/2},
 $$ 
we get that $ \|R-P\|\leq \eps$. Since $P$ is a projection and 
$0\leq R\leq 1$, we have that 
 $$
 \Spec(R)\subset [0,\eps]\cap [1-\eps,1].
 $$
If $Q$ denotes the spectral projection of $R$ corresponding to 
$[1-\eps,1]$, it follows that $Q\in \rho(D)'\cap B$ and 
 $$
 \|P-Q\|\leq \|P-R\|+\|R-Q\|\leq 2\eps.
 $$

Let $N$ be the dimension of $D$ and let $(e_i)_{i=1}^N$ be the 
family of minimal projections in $D$. Then the above $R$ is also 
given as $R=\sum_{i=1}^N \rho(e_i)P\rho(e_i)$. We set 
$a_i=\pi\alpha_f(e_i)$. Since 
 $$
 PR^nP=\sum_{i=1}^N a_i^{n+1}
 $$
for any $n\in\N$ and $n=0$, we get that 
 $$
 Pe^{itR}P=\sum_{i=1}^N\sum_{n=0}^\infty \frac{(it)^n}{n!} a_i^{n+1}
 =\sum_{i=1}^N a_ie^{ita_i}.
 $$

Since 
 $$
 \delta_\beta(P\rho(u)P\rho(u)^*P)=\delta_\beta(\pi(\alpha_f(u)\alpha_f(u)^*)) 
 =\pi(\delta_\alpha\alpha_f(u)\alpha_f(u)^*+\alpha_f(u)\delta_\alpha\alpha_f(u)^*)
 $$ 
for $u\in \U(D)$, we get that $\|\delta_\beta(PRP)\|< 2\eps$.
Similarly we get that 
 $$
 \|\delta_\beta(PR^nP)\|<(n+1)\eps.
 $$
This implies that $Pe^{itR}P$ belongs to the domain of 
$\delta_\beta$ for $t\in\R$ and that 
$\|\delta_\beta(Pe^{itR}P)\|\leq \eps (e^{|t|}+e^{|t|}|t|-1)$, 
which grows too rapidly in $t$. Since 
$\|\delta_\beta(e^{ita_i})\|\leq |t|\|\delta_\beta(a_i)\|\leq \eps 
|t|$, we get, from the expression of $Pe^{itR}P$ in terms of 
$a_i$'s, that 
 $$
 \|\delta_\beta(Pe^{itR}P)\|\leq N\eps(1+|t|),
 $$
which depends on the dimension $N$ of $D$. We shall give another 
estimate of $\|\delta_\beta(Pe^{itR}P)\|$, which grows 
polynomially in $t$ and is independent of $N$: 

\begin{lem} 
With $P,R$ as above and $t\in\R$ it follows that 
$\|\delta_\beta(Pe^{itR}P)\|\leq \eps(2|t|+t^2/2)$. 
\end{lem}
\begin{pf}
Note that 
 $$
 \delta_\beta(Pe^{itR}P)=
 \sum_{i=1}^N\sum_{n=0}^\infty\sum_{k=0}^n\frac{(it)^n}{n!}a_i^k\delta_\beta(a_i)a_i^{n-k}.
 $$
By using 
 $$
 \frac{1}{k!\cdot (n-k)!}\int_{0}^1s^k(1-s)^{n-k}ds=\frac{1}{(n+1)!},
 $$
the above equation equals 
 $$
 \sum_{i=1}^N\sum_{n=0}^\infty \sum_{k=0}^n \int_0^1ds \frac{(ist)^k(i(1-s)t)^{n-k}}{k!(n-k)!}
 (n+1)a_i^k\delta_\beta(a_i)a_i^{n-k}.
 $$
Substituting $n+1=1+k+\ell$ in the above formula, we get that 
 \BE
 \delta_\beta(Pe^{itR}P) &=&
    \sum_{k=0}^\infty\sum_{\ell=0}^\infty\sum_{i=1}^N \int_0^1ds \frac{(ist a_i)^k}{k!}
         \delta_\beta(a_i)\frac{(i(1-s)t a_i)^\ell}{\ell!} \\
  && +\sum_{k=1}^\infty\sum_{\ell=0}^\infty \sum_{i=1}^N \int_0^1 ds
          \frac{(ist a_i)^k}{(k-1)!}\delta_\beta(a_i)\frac{(i(1-s)t a_i)^{\ell}}{\ell!} \\
  && +\sum_{k=0}^\infty \sum_{\ell=1}^\infty\sum_{i=1}^N \int_0^1ds
          \frac{(ist a_i)^k}{k!}\delta_\beta(a_i)\frac{(i(1-s)t a_i)^\ell}{(\ell-1)!}.
 \EE       
Note that the sum over $i$ for $k=0$ and $\ell=0$ of the first 
term is zero because $\delta_\beta(\sum_ia_i)=\delta_\beta(P)=0$. 
We will evaluate the norm of the rest of the first term by 
splitting it into three terms $\Sigma_1=\sum_{k\geq 1,\ell\geq 
1}$, $\Sigma_3=\sum_{k=0,\ell\geq1}$, and $\sum_{k\geq 1,\ell=0}$, 
where $\sum_i\cdots$ is omitted. Similarly we will split the 
second term into two terms $\Sigma_2=\sum_{k\geq1,\ell\geq1}$, and 
$\sum_{k\geq1,\ell=0}$ and the third term into two 
$\sum_{k\geq1,\ell\geq1}\sum_i$ and 
$\Sigma_4=\sum_{k=0,\ell\geq1}\sum_i$ before evaluation. We can 
easily see that the unnamed terms can be expressed in terms of the 
named $\Sigma_1,\ldots,\Sigma_4$. 

We set 
 $$
 T=\sum_{i=1}^N \rho(e_i)P\delta_\beta(a_i)P\rho(e_i).
 $$
Since $\|\delta_\beta(a_i)\|\leq \eps$, we get that $\|T\|\leq 
\eps$. If $k\geq1$ and $\ell\geq 1$, it follows that 
 $$
 \sum_{i=1}^Na_i^k\delta_\beta(a_i)a_i^\ell=PR^{k-1}TR^{\ell-1}P.
 $$

We set 
 $$
 S=\sum_{i=1}^N P\delta_\beta(a_i)P\rho(e_i)=\sum_{i=1}^N\delta_\beta(a_i)P\rho(e_i).
 $$
If $\ell\geq 1$, it follows that 
 $$
 \sum_{i=1}^N \delta_\beta(a_i)a_i^\ell=SR^{\ell-1}P.
 $$                         
 
We assert that $\|S\|\leq \eps$. Since $SS^*=\sum_{i=1}^N 
\delta_\beta(a_i)a_i\delta_\beta(a_i)$ and $0\leq a_i\leq 1$, we 
get that 
 $$
 SS^*\leq \sum_{i=1}^N\delta_\beta(a_i)^2=XX^*,
 $$
where $X =(\delta_\beta(a_1),\delta_\beta(a_2),\ldots, 
\delta_\beta(a_N))\in M_{1N}(A)$. We naturally extend $\alpha$ to 
a flow on $M_{1 N}(A)\subset M_N(A)$, which we will also denote by 
$\alpha$. Similarly we extend $\pi$ to a representation of 
$M_N(A)=M_N\otimes A$ on $\C^N\otimes \Hil_\pi$. Then, since 
$\delta_\beta(a_i)=\pi\delta_\alpha\alpha_f(e_i)$, it follows that 
 $$
 X=\pi\delta_\alpha\alpha_f(e),
 $$
where $e=(e_1,e_2,\ldots,e_N)$, which has norm $1$. Since 
$\|\delta_\alpha\alpha_f\|\leq \eps$, which depends only on $f$, 
we get that $\|X\|\leq \eps$. This implies that $\|S\|\leq \eps$. 

Define a $C^\infty$-function $h$ on $\R$ by 
 $$
 h(t)=\frac{e^{it}-1}{it}=\sum_{k=1}^\infty\frac{(it)^{k-1}}{k!}
 $$                                  
and note that $|h(t)|\leq 1$. Since 
 $$
 \sum_{k=1}^\infty \frac{(ist)^k}{k!}R^{k-1}=ist \cdot h(stR),
 $$
we get that 
 \BE
  \Sigma_1(t)&\equiv&\sum_{k=1}^\infty\sum_{\ell=1}^\infty\sum_{i=1}^N 
  \int_0^1ds \frac{(ist a_i)^k}{k!}
         \delta_\beta(a_i)\frac{(i(1-s)t a_i)^\ell}{\ell!} \\
  &&= \int_0^1ds (-s(1-s)t^2)Ph(stR)Th((1-s)tR)P,
 \EE
which has the estimate $\|\Sigma_1(t)\|\leq \eps t^2/6$. We also 
have the following: 
 \BE
 \Sigma_2(t)&\equiv& \sum_{k=1}^\infty\sum_{\ell=1}^\infty\sum_{i=1}^N\int_0^1ds
 \frac{(ist a_i)^k}{(k-1)!}\delta_\beta(a_i)\frac{(i(1-s)t 
 a_i)^\ell}{\ell!} \\
 &=&\int_0^1ds (-s(1-s)t^2)Pe^{istR}Th((1-s)tR)P,
 \EE                                             
which has the estimate $\|\Sigma_2(t)\|\leq \eps t^2/6$. 
Furthermore we compute: 
 \BE
 \Sigma_3(t)&\equiv& \sum_{\ell=1}^\infty\sum_{i=1}^N
 \int_0^1ds \delta_\beta(a_i)\frac{(i(1-s)t a_i)^\ell}{\ell!}\\
 &=& \int_0^1ds \cdot i(1-s)t S h((1-s)tR)P,
 \EE
which implies that $\|\Sigma_3(t)\|\leq \eps |t|/2$. We also 
compute: 
 \BE
 \Sigma_4(t)&\equiv& \sum_{\ell=1}^\infty \sum_{i=1}^N \int_0^1ds 
    \delta_\beta(a_i)\frac{(i(1-s)t a_i)^\ell}{(\ell-1)!}\\
  &=& \int_0^1ds \cdot i(1-s)t Se^{i(1-s)tR}P,
 \EE
which implies that $\|\Sigma_4(t)\|\leq \eps |t|/2$. Since 
 $$
 \delta_\beta(Pe^{itR}P)=\Sigma_1(t)+\Sigma_2(t)+\Sigma_2(-t)^*
 +\Sigma_3(t)+\Sigma_3(-t)^*+\Sigma_4(t)+\Sigma_4(-t)^*,
 $$
we get the desired estimate. 
\end{pf}

Let $g$ be a $C^\infty$-function on $\R$ of compact support such 
that $g(t)=0$ for $t\in [0,1/3]$ and $g(t)=1$ for $t\in [2/3,1]$. 
We assume that $\eps$ is sufficiently small; in particular, 
$\eps<1/3$. Then $Q=g(R)$ is the spectral projection of $R$ 
corresponding to $[1-\eps,1]$. We have already noted that 
$\|P-Q\|\leq 2\eps$, which implies that $\Spec(PQP)\subset 
[1-2\eps,1]$ on $P\Hil_\rho=\Hil_\pi$. Thus the support projection 
of $PQP$ (resp. $QPQ$) is $P$ (resp. $Q$). 
 
We define a function $\hat{g}$ on $\R$ by 
 $$
 \hat{g}(s)=\frac{1}{2\pi}\int e^{-ist}g(t)dt.
 $$
Then we know that $\hat{g}$ is a rapidly decreasing 
$C^\infty$-function and that 
 $$
 g(R)=\int \hat{g}(t)e^{itR} dt.
 $$
Since $Pg(R)P\in D(\delta_\beta)$ and 
 $$
 \delta_\beta(Pg(R)P)=\int \hat{g}(t)\delta_\beta(Pe^{itR}P)dt,
 $$
we get, by the previous lemma, the estimate 
 $$
 \|\delta_\beta(Pg(R)P)\|\leq \eps\int |\hat{g}(t)|(2|t|+t^2/2)dt.
 $$

Suppose that we give such a function $g$ beforehand and let 
 $$
 C_1=\int |\hat{g}(t)|(2|t|+t^2/2)dt.
 $$
Then we choose $\eps>0$ so small that $C_1\eps$ is sufficiently 
small. Since $\Spec(PQP)\subset [2/3,1]$ (on 
$P\Hil_\rho=\Hil_\pi$), we have that $(PQP)^{-1/2}\in 
D(\delta_\beta)$. In the same way as above we have a constant 
$C_2$ such that 
 $$
 \|\delta_\beta((PQP)^{-1/2})\|\leq C_2\|\delta_\beta(PQP)\|\leq C_1C_2\eps,
 $$
which we assume is sufficiently small. 

We define a CP map $\Phi$ of $D$ into $PBP\subset \pi(A)$ by 
 $$
 \Phi(x)=(PQP)^{-1/2}PQ\rho(x)QP(PQP)^{-1/2},\ \ x\in D.
 $$  
Let $W=(PQP)^{-1/2}PQ\in B$, which is a partial isometry such that 
$WW^*=P$, $W^*W=Q$, and $\Phi(x)=W\rho(x)W^*,\ x\in D$. Since 
$[\rho(x),Q]=0$ for $x\in D$, $\Phi$ is a unital homomorphism. 
Since $\|PQ-P\|\leq 2\eps$, we have the estimate $\|W-P\|\leq 
3\eps$ (see 2.7 of \cite{Ch0}). Since $P\rho(x)P=\pi\alpha_f(x) 
\approx \pi(x)$ (up to $\eps^2/2$) we get for $x\in D$ with 
$\|x\|\leq 1$ that 
 $$
 \|\Phi(x)-\pi(x)\|\leq 3\eps+\eps^2/2.
 $$
Thus $\Phi$ is an injective homomorphism from $D$ into $PBP\subset 
\pi(A)$ such that $\Phi$ is close to $\pi|D$. 

Let $x=\sum_{i=1}^Nx_ie_i\in D$ be such that $\|x\|\leq 1$, i.e., 
$\max_i|x_i|\leq1$. For the same reasoning as for $Pe^{itR}P$ in 
the above lemma, we can conclude that $Pe^{itR}\rho(x)P\in 
D(\delta_\beta)$ with the same estimate 
 $$
 \|\delta_\beta(Pe^{itR}\rho(x)P)\|\leq \eps(2|t|+t^2/2).
 $$
(In the proof of the above lemma we just have to replace $T$ and 
$S$ by 
 $$ T'=\sum_{i=1}^N x_ie_iP\delta_\beta(a_i)P\rho(e_i)\ {\rm  and}\  
 S'=\sum_{i=1}^N x_iP\delta_\beta(a_i)P\rho(e_i),
 $$
respectively. Both $T'$ and $S'$ have the same estimates 
$\|T'\|\leq \eps$ and $\|S'\|\leq \eps$ as before since 
$\|x\|\leq1$.) Thus we get that $PQ\rho(x)QP\in D(\delta_\beta)$ 
and hence $\Phi(x)\in D(\delta_\beta)$; moreover, since 
$\|\delta_\beta((PQP)^{-1/2})\|\leq C_1C_2\eps$, 
$\|(PQP)^{-1/2}\|\leq (1-2\eps)^{-1/2}$, and 
$\|\delta_\beta(PQ\rho(x)P)\|\leq C_1\eps$, we have that  
 $$
 \|\delta_\beta(\Phi(x))\|\leq 2C_1C_2\eps(1-2\eps)^{-1/2}+(1-2\eps)^{-1}C_1\eps
 $$
for $x\in D$ with $\|x\|\leq 1$. 

Identifying $A$ with $\pi(A)$ and summing up the above arguments, 
we have proved the following assertion. For any $\eps>0$ there is 
an injective homomorphism $\Phi$ of $D$ into $A$ such that 
$\|\Phi-\id\|<\eps$ and $\Phi(D)\subset D(\delta_\alpha)$ and 
$\|\delta_\alpha|\Phi(D)\|<\eps$.  By 4.2 of \cite{Ch1} such a 
$\Phi$ is implemented by a unitary $w\in A$ such that $\|w-1\|\leq 
2\eps$ (and $\Phi(x)=wxw^*,\ x\in D$). 

Let $D_1=\Phi(D)$. We define $h=h^*\in A$ by 
 $$
 ih=\int_{\U(D_1)}\delta_\alpha(u)u^*d\nu(u), 
 $$ 
where $\nu$ is normalized Haar measure on the compact unitary 
group $\U(D_1)$. Then it follows that $\|h\|\leq \eps$ and 
$[ih,x]=\delta_\alpha(x),\ \ x\in D_1$. Let $u$ denote the 
differentiable $\alpha$ cocycle such that $ 
\frac{d}{dt}u_t|_{t=0}=-ih$. We set $v_t=w^*u_t\alpha_t(w)$, which 
is an $\alpha$-cocycle such that 
 $$
 \max_{|t|\leq 1}\|v_t-1\|\leq 2\|w-1\|+\eps\leq 5\eps
 $$
and 
 $$
 \Ad\,v_t\alpha_t(x)=\Ad\,w^*\Ad\,u_t\alpha_t\Phi(x)=\Ad\,w^*\Phi(x)=x,
 $$
for $x\in D$. This concludes the proof of Lemma \ref{I1}.  
\end{pf}

It is instructive to see what $\Phi$ is. By computation, since 
$Q=g(R)$, we have that $PQ\rho(e_i)P=a_ig(a_i)$. Since 
$\Spec(a_i)\subset [0,\eps]\cup[1-\eps,1]$, if $g_1$ is a 
continuous function on $\R$ such that $g=0$ on $[0,1/3]$ and 
$g(t)=t$ on $[2/3,1]$, we get $g_1(a_i)=a_ig(a_i)$. Thus it 
follows, with such a $g_1$, that 
 $$
 \Phi(e_i)=b_{g_1}^{-1/2}g_1(\alpha_f(e_i))b_{g_1}^{-1/2},
 $$
where $b_{g_1}=\sum_{j=1}^Ng_1(\alpha_f(e_j))$. (Since there is  
freedom in the above proof, we may as well take $g$ instead of 
$g_1$ in the above formula (with a small but different 
$\|\delta_\alpha|\Phi(D)\|$). Thus the spectral projections 
$q_i$'s of $\alpha_f(e_i)$'s corresponding to $[1-\eps,1]$ may not 
be mutually orthogonal, but 
$(\sum_jq_j)^{-1/2}q_i(\sum_jq_j)^{-1/2}$'s are mutually 
orthogonal projections.)  

\medskip
\noindent {\em Proof of Proposition \ref{Inv}}

What we have to show is that Lemma \ref{I1} follows without the 
assumption that $D$ is abelian. Hence suppose that $D$ is a unital 
finite-dimensinal C$^*$-subalgebra of $A$ such that 
 $$\sup_{|t|\leq 1}\|(\alpha_t-\id)|D\|<\delta,
 $$
where $\delta$ is the one obtained for $\eps$ in Lemma \ref{I1}. 
Let $Z$ be a maximal abelan C$^*$-subalgebra of $D$. Then by Lemma 
\ref{I1} we get an $\alpha$-cocycle $u$ such that 
$\Ad\,u_t\alpha_t|Z=\id$ and $\max_{|t|\leq1}\|u_t-1\|<\eps$. 

Let $(e_{ij}^{(k)})$ be a family of matrix units of $D$ such that 
the linear span of $e_{ii}^{(k)}$ for all $i$ and $k$ equals $Z$. 
We set 
 $$
 v_t=\sum_k\sum_ie_{i1}^{(k)}\Ad\,u_t\alpha_t(e_{1i}^{(k)}),
 $$ 
for $t\in\R$,  which is an $\Ad\,u\alpha$-cocycle. Then we have 
that 
 $$
 \Ad(v_t 
 u_t)\alpha_t(e_{ab}^{(c)})=e_{a1}^{(c)}\Ad\,u_t\alpha_t(e_{1a}^{(c)} 
 e_{ab}^{(c)}e_{b1}^{(c)})e_{1b}^{(c)}=e_{ab}^{(c)},
 $$
i.e., $\Ad(v_tu_t)\alpha_t|D=\id$. Since $\|v_t-1\|$ is equal to 
 $$
 \|\sum_{i,k}e_{i1}^{(k)}(\Ad\,u_t\alpha_t(e_{1i}^{(k)})-e_{1i}^{(k)})\|
 =\max_{i,k}\|\Ad\,u_t\alpha_t(e_{1i}^{(k)})-e_{1i}^{(k)}\|\leq
  \|(\Ad\,u_t\alpha_t-\id)|D\|,
 $$
it follows that $ \max_{|t|\leq 1}\|v_tu_t-1\|\leq 3\eps+\delta$.
Since $t\mapsto v_tu_t$ is an $\alpha$-cocycle, this concludes the 
proof. 

\medskip  

Now we turn to the proof of Lemma \ref{M1}. We start with the 
following lemma. 

\begin{lem}   \label{M2}
Let $\alpha$ be a flow on a unital AF algebra $A$. Suppose that 
there is an increasing sequence $(A_n)$ of finite-dimensional  
C$^*$-subalgebras of $A$ such that $\bigcup_nA_n$ is dense in $A$ 
and 
 $$
 \sup_{t\in [0,1]}\dist(\alpha_t(A_n),A_n)\ra0.
 $$                                            
Then for any $\eps>0$ there is an $\alpha$-cocycle $u$ and a 
subsequence $(n_i)$ in $\N$ such that 
$\Ad\,u_t\alpha_t|A_{n_i}\cap A_{n_i}'=\id$ for all $i$ and
$\max_{t\in [0,1]}\|u_t-1\|<\eps$. 
\end{lem}                                         
\begin{pf}
Let $e$ be a non-zero central projection of $A_n$. Then for any 
projection $p\in A_n$ different from $e$, we have that 
$\|e-p\|=1$. 

If $\dist(\alpha_t(A_n),A_n)<\delta$ for $t\in [0,1]$ for a small 
$\delta$, then it follows that $\|\alpha_t(e)-e_t\|<\delta$ for 
some $e_t\in A_n$. We may assume that $e_t$ is a projection by 
replacing $\delta$ by $2\delta$. Then from the above remark we get 
that $e_t=e$. Thus we have that $\|\alpha_t(e)-e\|<2\delta$ for 
all projections $e\in A_n\cap A_n'$ and for $t\in[0,1]$. Since any 
element $x\in A_n\cap A_n'$ with $0\leq x\leq1$ is a convex 
combination of projections in $A_n\cap A_n'$, it follows that 
$\sup_{t\in [0,1]}\|(\alpha_t-\id)|A_n\cap A_n'\|$ converges to 
zero as $n\ra\infty$. 
                          
By Lemma \ref{I1} there is an $n_1\in \N$ and an $\alpha$-cocycle 
$u^{(1)}$ such that $\max_{t\in[0,1]}\|u_t^{(1)}-1\|<\eps/2$ and 
$\Ad\,u^{(1)}_t\alpha_t|A_{n_1}\cap A_{n_1}'=\id$. Let 
$Z_1=A_{n_1}\cap A_{n_1}'$. Since $\Ad\,u^{(1)}\alpha$ leaves 
$A\cap Z_1'$ invariant, $\bigcup_{n>n_1} A_n\cap Z_1'$ is dense in 
$A\cap Z_1'$, and 
$\sup_{t\in[0,1]}\dist(\Ad\,u_t^{(1)}\alpha_t(A_n\cap 
Z_1'),A_n\cap Z_1')\ra0$ as $n\ra\infty$, we can repeat this 
argument for $\Ad\,u^{(1)}\alpha|A\cap Z_1'$ and $(A_n\cap Z_1' 
)_{n>n_1}$ with $\eps$ replaced by $\eps/2$. Thus we have an 
$n_2>n_1$ and an $\Ad\,u^{(1)}\alpha$-cocycle $u^{(2)}$ in $A\cap 
Z_1'$ such that $\Ad(u^{(2)}_tu^{(1)}_t)\alpha_t$ is the identity 
on the center of $A_{n_2}\cap Z_1'$ and 
$\max_{t\in[0,1]}\|u_t^{(2)}-1\|<\eps/4$. It then follows that 
$\Ad(u_t^{(2)}u_t^{(1)})\alpha_t|A_{n_i}\cap A_{n_i}'=\id$ for 
$i=1,2$ and that $t\mapsto u_t^{(2)}u_t^{(1)}$ is an 
$\alpha$-cocycle such that 
$\max_{t\in[0,1]}\|u_t^{(2)}u_t^{(1)}-1\|<\eps/2+\eps/4$. In this 
way we can complete the proof. 
\end{pf}                                                           

Now we assume that $\alpha$ fixes each central projection of $A_n$ 
for all $n$. Let $C$ denote the C$^*$-subalgebra generated by 
$\bigcup_n (A_n\cap A_n')$, which is an AF abelian \cstar. Let 
$B=A\cap C'$. Since $\alpha$ fixes each element of $C$, $\alpha$ 
restricts to a flow on $B$, which we will denote by $\beta$. 

Let $B_n=A_n\cap C'$, which is the relative commutant of the 
C$^*$-subalgebra $C_n$ in $A_n$, where $C_n$ is generated by 
$\bigcup_{k=1}^n(A_k\cap A_k')$. Thus there is a norm one 
projection of $A$ onto $A\cap C'$, sending $A_n$ onto $B_n=A_n\cap 
C_n'$ and $B$ is an AF algebra. We identify $C$ with the 
continuous functions $C(\Gamma)$, where $\Gamma$ is the compact 
Hausdorff space of characters of $C$. Then we can regard $B$ as 
the \cstar\ of continuous sections over $\Gamma$; the fiber at 
$\gamma\in \Gamma$ will be denoted by $B^\gamma$, which is a UHF 
algebra (or a matrix algebra), and the canonical map of $B$ onto 
$B^\gamma$ will be denoted by $\Phi^\gamma:x\mapsto x(\gamma)$. To 
see what $B^\gamma$ is, we find a decreasing sequence $(e_n)$ of 
projections such that $e_n\in C_n$ is minimal with 
$\gamma(e_n)=1$; then $B^\gamma$ is obtained as the inductive 
limit of the sequence $e_1A_1e_1\ra e_2A_2e_2\ra 
e_3A_3e_3\ra\cdots$, where the map of $e_nA_ne_n=B_ne_n=(A_n\cap 
C_n')e_n$ into $e_{n+1}A_{n+1}e_{n+1}$ is given by $x\mapsto 
xe_{n+1}$. Let $B^\gamma_n$ denote the image of $e_nA_ne_n$ in 
$B^\gamma$. We define a flow $\beta^\gamma$ on $B^\gamma$ by the 
requirement that 
$\Phi^\gamma(\beta_t(x))=\beta^\gamma_t(\Phi^\gamma(x)),\ x\in B$.

Since $\sup_{t\in[0,1]}\dist(\alpha_t(A_n),A_n)\ra0$, we obtain 
that                                                    
 $$
 \sup_{t\in[0,1]}\dist(\beta^\gamma_t(B^\gamma_n),B^\gamma_n)\ra0
 $$
for any $\gamma\in\Gamma$. Thus, by \ref{M3}, each $\beta^\gamma$ 
is a cocycle perturbation of a UHF flow. To show that $\alpha$ is 
a cocycle perturbation of an AF flow, we would have to use the 
fact that the convergence in the above display is uniform in 
$\gamma\in\Gamma$.

We shall prove the first half of Lemma \ref{M1} as follows. 
Suppose, in particular, that 
$\sup_{t\in[0,1]}\dist(\alpha_t(A_1),A_1)$ is sufficiently small. 
We have to show that there is an $\alpha$-cocycle $w$ such that 
$\sup_{t\in[0,1]}\|w_t-1\|$ is arbitrarily small and 
$\Ad\,w_t\alpha_t(A_1)=A_1$. 

Suppose that the center of $A_1$ is $K$-dimensional, being spanned 
by minimal central projections $z_1,\ldots,z_K$, and choose 
$\gamma_1,\gamma_2,\ldots,\gamma_K\in\Gamma$ such that  
$\gamma_i(z_i)=1$ for each $i=1,2,\ldots,K$. We apply Lemma 
\ref{M3} to each $(B^{\gamma_i},\beta^{\gamma_i})$ (or apply 
Prop.~\ref{A} if $B^{\gamma_i}$ is a matrix algebra); thus we get 
a $\gamma_i$-cocycle $u^i$ such that $\|u^i_t-1\|$ is very small 
for $t\in[0,1]$ and $\Ad\,u^i_t\beta^{\gamma^i}_t$  fixes 
$B^\gamma_1$. What is important here is there is a pure ground 
state $\phi_i$ for $(B^{\gamma_i},\Ad\,u^i\beta^{\gamma_i})$, 
which is still pure on $B^\gamma_1$. We may regard $\phi_i$ as a 
pure state $\phi_i\Phi^{\gamma_i}$ on $z_iAz_i$, which is pure on 
$z_iA_1z_i(\supset B^\gamma_1)$. We consider the system 
$(z_iAz_i,\alpha^i=\alpha|z_iAz_i)$. By lifting $u^i$ (see below), 
we then find an $\alpha^i$-cocycle $v^i$ such that 
$\max_{t\in[0,1]}\|v^i_t-1\|$ is sufficiently small and $\phi_i$ 
is $\Ad\,v^i_t\alpha^i_t$-invariant. Since 
$\dist(\Ad\,v^i_t\alpha^i_t(A_1z_i),A_1z_i)\leq 
2\|v^i_t-1\|+\dist(\alpha_t(A_1),A_1)$, which is very small, we 
apply Lemma \ref{C} to get an $\alpha^i$-cocycle $w^i$ such that 
$\max_{t\in[0,1]}\|w_t-1\|$ is small and $\Ad\,w^i_t\alpha^i$ 
fixes $A_1z_i(\subset z_iAz_i$). By combining these $w^i$ we get 
an $\alpha$-cocycle $w$ such that $\max_{t\in[0,1]}\|w_t-1\|$ is 
small and $\Ad\,w_t\alpha_t$ fixes $A_1$.   
                   
\begin{lem}
Let $\gamma\in \Gamma$ and let $u$ be a $\beta^\gamma$-cocycle. 
Then there is an $\alpha$-cocycle $v$ such that $v_t\in B=A\cap 
C'$ and $\Phi^\gamma(v_t)=u_t,\ t\in\R$. If  $\sup_{t\in 
[0,1]}\|u_t-1\|<\delta$ holds,  $\sup_{t\in[0,1]}\|v_t-1\|<\delta$ 
can be imposed. 
\end{lem}
\begin{pf}
We find a $w\in \U(B^\gamma)$ and a differentiable 
$\beta^\gamma$-cocycle $z$ such that $\|w-1\|$ is small and 
$u_t=wz_t\beta^\gamma_t(w^*)$. We can then find a $W\in B=A\cap 
C'$ and $H\in B_{sa}$ such that $\Phi^\gamma(W)=w$ and 
$\Phi^\gamma(H)=-id z_t/dt|_{t=0}$. Then we find a $\beta$-cocycle 
$Z$ by solving the equation $ dZ_t/dt=Z_t\beta_t(iH)$ with 
$Z_0=1$, which satisfies that $\Phi^\gamma(Z_t)=z_t$. We set 
$v_t=WZ_t\beta_t(W^*)$, which is an $\alpha$-cocycle with 
$\Phi^\gamma(v_t)=u_t$. 

Since $\sigma\mapsto \|\Phi^\sigma(v_t)-1\|$ is continuous, 
$t\mapsto \|\Phi^\sigma(v_t)-1\|$ is equi-continuous in 
$\sigma\in\Gamma$, and $\Gamma$ is totally disconnected, the last 
condition is satisfied by replacing $v$ by $t\mapsto v_te+1-e$, 
where $e\in C$ is a projection with $\gamma(e)=1$. 
\end{pf}

Let $\eps_1>0$. For $\eps=\eps_1/2$ we choose a 
$\delta_1\equiv\delta>0$ as in Lemma \ref{C}, where we may assume 
that $\delta_1<\eps_1/2$.  For $\eps=\delta_1/3>0$ we choose a 
$\delta_2\equiv \delta>0$ as in Lemma \ref{M3} (or Prop.~\ref{A}). 
We may assume that $\delta_2<\delta_1$. 

Suppose that $\sup_{t\in[0,1]}\dist(\alpha_t(A_1),A_1)<\delta_2$. 
Let $z_1,z_2,\ldots,z_K$ be the minimal central projections in 
$A_1$ and let $\gamma_1,\ldots,\gamma_K\in \Gamma$ be such that 
$\gamma_i(z_i)=1$. Then we apply \ref{M3} to 
$(B^{\gamma_i},\beta^{\gamma_i})$ to get a 
$\beta^{\gamma_i}$-cocycle $u^i$ such that $\sup_{t\in 
[0,1]}\|u_t^i-1\|<\delta_1/3$ and 
$\Ad\,u^i_t\beta^{\gamma_i}_t(B^{\gamma_i}_1)=B^{\gamma_i}_1$. 
Since $\beta^{\gamma^i}$ is a UHF flow, there is a pure ground 
state $\phi_i$ for $\Ad\,u^i\beta^{\gamma_i}$. Note that 
$\phi_i|B^{\gamma_i}_1$ is pure. We regard $\phi_i$ as a state on 
$z_iAz_i$ by denoting $\phi_i\Phi^{\gamma_i}|z_iAz_i$ again by 
$\phi_i$. Note that $\phi_i|z_iA_1z_i$ is pure. We lift $u^i$ to 
an $\alpha$-cocycle $v^i$ such that $v^i_t\in B$ and 
$\sup_{t\in[0,1]}\|v^i_t-1\|<\delta_1/3$. We regard $v^i$ as a 
cocycle in $z_iAz_i$ with respect to $\alpha|z_iAz_i$. Then it 
follows that 
 $$
 \sup_{t\in[0,1]}\dist(\Ad\,v^i_t\alpha_t(z_iA_1),z_iA_1)
 <\delta_1.
 $$

Then applying \ref{C} to $(z_iAz_i, \Ad\,v^i_t\alpha_t|z_iAz_i)$ 
with a pure invariant state $\phi_i$, we get a cocycle $w^i$ with 
respect to $\Ad\,v_t^i\alpha_t|z_iAz_i$ such that 
 $$\sup_{t\in[0,1]}\|w_t-1\|<\eps/2 
 \ \ {\rm and}\ \  \Ad(w^i_t 
 v^i_t)\alpha_t(z_iA_1)=z_iA_1.
 $$
Note that $t\mapsto w^i_tv^i_t$ is an $\alpha$-cocycle in 
$z_iAz_i$ such that $\sup_{t\in [0,1]}\|w^i_tv^i_t-1\|<\eps_1$. 
Then $t\mapsto \sum_iw^i_tv^i_t$ is the desired $\alpha$-cocycle. 
This completes the proof of the first half of Lemma \ref{M1}. 

To prove the latter part, let $A_0$ be a C$^*$-subalgebra of $A_1$ 
such that $\alpha_t(A_0)=A_0$. We have to show that the above
$\alpha$-cocycle can be chosen from $A\cap A_0'$.

Since $\dist(\alpha_t(A_n\cap A_0'),A_n\cap A_0')\leq 
\dist(\alpha_t(A_n),A_n)$, we apply the first part to 
$\alpha|A\cap A_0'$ to find an $\alpha$-cocycle $u$ in $A\cap 
A_0'$ such that $\sup_{t\in[0,1]}\|u_t-1\|<\eps$ and 
$\Ad\,u_t\alpha_t(A_1\cap A_0')=A_1\cap A_0'$. Since 
$\Ad\,u_t\alpha_t(A_0)=A_0$, we get that $\Ad\,u_t\alpha_t(A_1\cap 
Z_0')=A_1\cap Z_0'$, where $Z_0=A_0\cap A_0'$. Note that 
$\sup_{t\in[0,1]}\dist(\Ad\,u_t\alpha_t(A_1),A_1)< 
\delta_1+2\eps$. 

Let $z_1,\ldots,z_K$ be the minimal central projections in $A_1$ 
as before. If $(A_1\cap Z_0')z_i$ is the direct sum of $L_i+1$ 
factors, there are $L_i$ partial isometries $w_1,\ldots,w_{L_i}$ 
in $A_1z_i$ such that $w_j^*w_j=w_1^*w_1$ is a minimal projection 
invariant under $\Ad\,u\alpha$, $w_jw_j^*$ is a minimal projection 
invariant under $\Ad\,u\alpha$, and $A_1z_i$ is generated by 
$(A_1\cap Z_0')z_i$ and $w_j,\ j=1,\ldots,L_i$.      

Since $\Ad\,u_t\alpha_t(w_j)$ is almost contained in $A_1$ (up to 
the order of $\delta_1+2\eps$) for $t\in[0,1]$, it follows that 
$\|\Ad\,u_t\alpha_t(w_j)-c_tw_j\|<\delta_1+2\eps$ for some 
$c_t\in\C$ for $t\in[0,1]$. Note that $t\mapsto 
w_j\Ad\,u_t\alpha_t(w_j^*)$ is a cocycle with respect to 
$\Ad\,u\alpha|w_jw_j^*Aw_jw_j^*$.       

\begin{lem}  \label{X}
For any $\eps>0$ there exists a $\delta>0$ satisfying the 
following condition: If $v$ is an $\alpha$-cocycle such that 
$v_t\stackrel{\delta}{\in}\C1$ for $t\in[0,1]$, then there is a 
$p\in\R$ such that $\|v_t-e^{ipt}1\|<\eps$ for $t\in [0,1]$. 
\end{lem}           

Thus we find a $\lambda_j\in\R$ such that 
$\sup_{t\in[0,1]}\|w_j\Ad\,u_t\alpha_t(w_j^*)-e^{i\lambda_jt}w_jw_j^*\|
\approx0$. We extend a cocycle $t\mapsto 
e^{-i\lambda_jt}w_j\Ad\,u_t\alpha_t(w_j^*)$ in $w_jw_j^*(A\cap 
(A_1\cap Z_0')')$ to a cocycle $v^j$ in $p_j(A\cap (A_1\cap 
Z_0')')$, where $p_j$ is the central support projection of 
$w_jw_j^*$ in $A_1\cap Z_0'$ (or in $z_i(A_1\cap Z_0')$). Note 
that 
$\|v^j_t-p_j\|=\|e^{-i\lambda_jt}w_j\Ad\,u_t\alpha_t(w_j^*)-w_jw_j^*\|$. 
We define $v_t=\sum_{j}v^j_t+p_0$, where $p_0$ is the central 
support of $w_j^*w_j=w_1^*w_1$ in $z_i(A_1\cap Z_0')$. Then $v$ is 
a cocycle in $z_iA\cap (A_1\cap Z_0')'$ with respect to 
$\Ad\,u\,\alpha$, $\sup_{t\in[0,1]}\|v_t-z_i\|\approx0$, and 
 $$
 \Ad(v_tu_t)\alpha_t(w_j)=e^{-i\lambda_jt}w_j.
 $$ 
Thus, since $A_1\cap Z_0'$ and all $w_j$ generate $z_iA_1$ and 
$A_1\cap Z_0'\supset A_0$, we get that 
 $$\Ad(v_tu_t)\alpha_t(z_iA_1)=z_iA_1
 \ \ {\rm and}\ \  
 \Ad(v_tu_t)\alpha_t(z_iA_0)=z_iA_0.
 $$
We apply this argument to each $z_i$. This concludes the proof of 
Lemma \ref{M1}.       

\medskip

\noindent 
 {\em Proof of \ref{X}}      
We suppose that for a small $\delta>0$ there is a continuous 
function $f:[0,1]\ra\R$ such that $\|u_t-e^{if(t)}1\|<\delta$ and 
$f(0)=0$. For $t_1,t_2\in[0,1]$ with $t_1+t_2\leq1$, we have that 
$\|u_{t_1}\alpha_{t_1}(u_{t_2})-e^{if(t_1)+if(t_2)}1\|<2\delta$, 
which implies that $|e^{i(f(t_1+t_2)-f(t_1)-f(t_2))}-1|<3\delta$.
Letting $\delta_0=\arcsin 3\delta$, we get
 $$
 |f(t_1+t_2)-f(t_1)-f(t_2)|<\delta_0.
 $$
We replace $f$ by $f(t)-f(1)t$, which still satisfies the above 
inequality for $t_1,t_2,t_1+t_2\in[0,1]$. With this replacement we 
have assumed that $f(0)=0=f(1)$. 

Let $\mu=\max\{|f(t)|\ |\ t\in [0,1]\}$. Suppose that 
$\mu>3\delta_0$. There is an $s\in [0,1]$ such that $|f(s)|=\mu$. 
If $s\leq 1/2$, then 
$|f(2s)|>2|f(s)|-\delta_0=2\mu-\delta_0>\mu+2\delta_0$, which is a 
contradiction. Hence we have that $s>1/2$. Then 
$|f(1-s)|>|f(s)|-\delta_0$, and hence 
$|f(2(1-s))|>2|f(1-s)|-\delta_0>2|f(s)|-3\delta_0>\mu$, which is 
again a contradiction. Thus we should have that $\mu\leq 
3\delta_0$. This implies that
 $$
 \|u_t-e^{if(1)t}1\|<|e^{if(t)}-e^{if(1)t}|+\delta<3\delta_0+\delta.
 $$
Since $\delta_0\approx 3\delta$, this concludes the proof.

\medskip
\small

\end{document}